\documentclass[11pt]{amsart}
\usepackage[top=2cm,bottom=2cm,left=1.5cm,right=1.5cm,marginparwidth=1.75cm]{geometry}

\usepackage[square,compress,comma, numbers,sort]{natbib}
\usepackage[colorlinks=true, citecolor=red, linkcolor=blue]{hyperref}
\usepackage{amsfonts,mathtools,mathabx,amssymb}
\usepackage{color}

 \usepackage{xparse}

\usepackage[shortlabels]{enumitem}

\usepackage{amsthm}
\usepackage{cleveref}
\usepackage{aliascnt}

\newtheorem{theorem}{Theorem}[section]

\newaliascnt{lemma}{theorem}
\newtheorem{lemma}[lemma]{Lemma}
\aliascntresetthe{lemma}
\crefname{lemma}{Lemma}{Lemmas}
\Crefname{lemma}{Lemma}{Lemmas}

\newaliascnt{proposition}{theorem}
\newtheorem{proposition}[proposition]{Proposition}
\aliascntresetthe{proposition}
\crefname{proposition}{Proposition}{Propositions}
\Crefname{proposition}{Proposition}{Propositions}

\newaliascnt{corollary}{theorem}
\newtheorem{corollary}[corollary]{Corollary}
\aliascntresetthe{corollary}
\crefname{corollary}{Corollary}{Corollaries}
\Crefname{corollary}{Corollary}{Corollaries}

\theoremstyle{definition}
\newaliascnt{definition}{theorem}

\aliascntresetthe{definition}
\crefname{definition}{Definition}{Definitions}
\Crefname{definition}{Definition}{Definitions}

\newaliascnt{example}{theorem}
\newtheorem{example}[example]{Example}
\aliascntresetthe{example}
\crefname{example}{Example}{Examples}
\Crefname{example}{Example}{Examples}

\theoremstyle{remark}
\newaliascnt{remark}{theorem}
\newtheorem{remark}[remark]{Remark}
\aliascntresetthe{remark}
\crefname{remark}{Remark}{Remarks}
\Crefname{remark}{Remark}{Remarks}

\allowdisplaybreaks[4]

\DeclareMathOperator*{\cov}{cov}

\DeclareMathOperator*{\Cov}{Cov}
\DeclareMathOperator*{\Var}{Var}

\DeclarePairedDelimiterXPP\pk[1]{\mathbb{P}}\{ \}{}{ #1}
\DeclarePairedDelimiterXPP\E[1]{\mathbb{E}}\{ \}{}{	#1}

\NewDocumentCommand{\ceil}{s O{} m}{%
  \IfBooleanTF{#1} 
    {\left\lceil#3\right\rceil} 
    {#2\lceil#3#2\rceil} 
}
\NewDocumentCommand{\floor}{s O{} m}{%
  \IfBooleanTF{#1} 
    {\left\lfloor#3\right\rfloor}
    {#2\lfloor#3#2\rfloor}
}

\definecolor{c20}{rgb}{0.,0.7,0.}
\definecolor{c30}{rgb}{0.,0.,1.}
\definecolor{c40}{rgb}{1,0.1,0.7}
\definecolor{c50}{rgb}{1,0,0}
\definecolor{c60}{rgb}{1,0.9,0.1}
\definecolor{c70}{rgb}{0.50,1.00,0.00}

\def\eH#1{{\textcolor{c30}{#1}}}
\def\eH#1{#1}

\def\eH#1{#1}

\definecolor{orange}{rgb}{1,0.5,0}
\def\k2#1{{\textcolor{orange}{#1}}}

\def\aa#1{{\textcolor{black}{#1}}}
\def\xx#1{{\textcolor{black}{#1}}}

\def\ks#1{{\textcolor{black}{#1}}}

\def\kE#1{{\textcolor{black}{#1}}}

\newcommand{\BS}{\begin{proposition}}
\newcommand{\ES}{\end{proposition}}
\newcommand{\BT}{\begin{theorem}}
\newcommand{\ET}{\end{theorem}}
\newcommand{\BK}{\begin{corollary}}
\newcommand{\EK}{\end{corollary}}

\newcommand{\BEX}{\begin{example}}
\newcommand{\EEX}{\end{example}}

\newcommand{\BRM}{\begin{remark}}
\newcommand{\ERM}{\end{remark}}

\newcommand{\BEL}{\begin{lemma}}
\newcommand{\EEL}{\end{lemma}}

\newcommand{\prooftheo}[1]{Proof of  Theorem \ref{#1}:}
\newcommand{\proofprop}[1]{Proof of  Proposition \ref{#1}:}
\newcommand{\prooflem}[1]{Proof of  Lemma \ref{#1}:}
\newcommand{\proofkorr}[1]{Proof of  Corollary \ref{#1}:}

\newcommand{\QED}{\hfill $\Box$}

\newcommand{\COM}[1]{}

\def\IF{\infty}

\newcommand{\R}{\mathbb{R}}
\newcommand{\inr}{\in \R}

\newcommand{\BQN}{\begin{eqnarray}}
\newcommand{\EQN}{\end{eqnarray}}
\newcommand{\BQNY}{\begin{eqnarray*}}
\newcommand{\EQNY}{\end{eqnarray*}}

\newcommand{\limit}[1]{\lim_{#1 \to   \infty}}

\def\bqny#1{\begin{eqnarray*} #1 \end{eqnarray*}}
\def\bqn#1{\begin{eqnarray} #1 \end{eqnarray}}

\newcommand\ind[1]{\mathbb{I}{\left\{#1\right\}}}

\def\kTt{\kappa_V(t,T)}

\begin{document}
\title
[  Expected Infimum \& persistence probabilities for stationary BR-Processes]
{  Expected Infimum and persistence probabilities of Log-Normal  Stationary Brown-Resnick Processes}

\author{Krzysztof D\c{e}bicki}
\address{Krzysztof D\c{e}bicki, Mathematical Institute, University of Wroc\l aw, pl. Grunwaldzki 2/4, 50-384 Wroc\l aw, Poland}
\email{Krzysztof.Debicki@math.uni.wroc.pl\ \ {\it{(This is the corresponding author)}}}

\author{Enkelejd  Hashorva}
	\address{Enkelejd Hashorva, Department of Actuarial Science,
		University of Lausanne,
		UNIL-Dorigny, 1015 Lausanne, Switzerland
	}
	\email{Enkelejd.Hashorva@unil.ch}

\author{Svyatoslav Novikov}
	\address{Svyatoslav Novikov, Department of Actuarial Science,
		University of Lausanne,
		UNIL-Dorigny, 1015 Lausanne, Switzerland
	}
	\email{Svyatoslav.Novikov@unil.ch}

\bigskip

\date{\today}
 \maketitle

\begin{quote}

We investigate the asymptotics of the expected infimum of log-normal Brown-Resnick stationary processes,
a class of processes that arise naturally in the study of extremes of Gaussian processes and max-stable processes.
Specifically, we analyse the functional
$$\mathcal{G}_V(T) = \mathbb{E}\left\{\inf_{t \in [0,T]}e^{ \sqrt{2}V(t)-\sigma^2_V(t)} \right\}, \quad T>0,
$$
where $V$ is  a  centered Gaussian process
with stationary increments, continuous sample paths and variance $\sigma_V^2$,
\aa{and a closely related problem of} the
decay rate of the persistence probability
$$p_V(T,C)=\pk*{\inf_{t\in [0,T]} (\sqrt{2} V(t)- \sigma^2_V(t)) > C}$$
for some constant $C<0$.\\
  For both $\mathcal{G}_V(T)$ and $p_V(T,C)$ we derive exact asymptotics as $T \to \infty$
  for a broad class of processes $V$, including fractional Brownian motion
  with Hurst parameter $H \in (1/2,1]$.
  \aa{For the latter, the behavior in the short-range dependence regime $H \in (0, 1/2]$ is markedly different and, in general, more delicate;
  we find logarithmic asymptotics for a family of processes that includes this case.}
   Our results provide sharp bounds and comparison principles,
   and highlight the contrast between the behavior of infimum and
   supremum functionals for such processes.
   The  discrete-time analogues and connections to Pickands constants
   are also discussed.\\
\end{quote}

{\it Key words}: Brown-Resnick stationarity;
infimum functional; persistence probabilities;
Gaussian processes; fractional Brownian motion;
Pickands constant;  Berman condition; stationary increments. \\

{\it MSC Subject Classification}:
60G15; 60G22; 60G70; 60G18; 60F05; 60F10; 60F15; 60F17; 60F99.

\section{Introduction}
Let $V(t),t\inr $ be a centered,   real-valued Gaussian process with stationary increments,   variance function $\sigma^2_V(t),t\inr$ and continuous sample paths.
Following \cite{kab2009}, we
  call the  log-normal process
\[
Z_V(t)=\exp(\sqrt{2}V(t)-\sigma^2_V(t)),\qquad t\in \R
\]
{\it Brown-Resnick stationary}. As shown therein,   the corresponding max-stable $X(t), t\inr$ defined
via its de Haan representation   (see e.g., \cite{kulik:soulier:2020,Hrovje}) by
\bqn{\label{eq1}
	{X} (t)=  \max_{i\ge 1} \Bigl(\sum_{k=1}^i \mathcal{E}_k\Bigr)^{-1}  Z_V^{(i)}(t), \quad t \inr
}
is stationary. Here   $\mathcal{E}_k, k\ge 1$ are unit iid exponential random variables (rv's) being independent of any other random element and $Z_V^{(i)}$'s are iid copies of $Z_V$.
A canonical class of examples of $X$ is for $V=B_H$, with $B_H(t),t\inr $ a standard fractional Brownian motion (fBM) with Hurst parameter $H\in (0,1]$.

The family of {\it Brown-Resnick stationary} processes, initially  introduced in \cite{bro1977,eddy1980distribution}
for $V=B_{1/2}$  and $V=B_1$, respectively,  plays an important  role in the asymptotic theory of heavy-tailed time series,
see e.g.,   \cite{kulik:soulier:2020, Resnickart, Ilya25}. Moreover, it has  been instrumental in deriving
novel representations of extremal indices including Pickands constants.  More specifically, the functional
\begin{eqnarray}
	\mathcal{H}_V([0,T])=\E*{\sup_{t\in [0,T]}Z_V(t)}, \quad  T \in [0, \IF) \label{wills}
\end{eqnarray}
has been studied extensively  in the literature
in various contexts including  Pickands constants,   and extremal indices, see e.g.,
 \cite{MR3745388,debicki2017approximation,ZKE} and the references therein. \\

As shown in \cite{Htilt},   the subadditivity property of supremum functional and Borell-TIS inequality imply
\bqn{ \limit{T} \frac{\mathcal{H}_V([0,T])}{T}= \inf_{T>0} \frac{\mathcal{H}_V([0,T])}{T} = \mathcal{H}_V<\IF.
\label{Ph}
}
The Pickands constant $\mathcal{H}_V$ is positive if
\[
Z_V(t)\to0\quad\text{almost surely as }|t|\to\infty.
\]
Under the present assumptions, this condition characterises pure
dissipativity of the max-stable process $X$; see, e.g.,
\cite{Genna04,Genna04c,WangStoev,kab2009,Hrovje,kulik:soulier:2020,hashorva2021shiftinvariant}.
Moreover,  as advocated in \cite{kab2009},
$\mathcal{H}_V$ depends only on
the variogram $\gamma_V$ defined by
 \[
\gamma_V(h) =
\frac{1}{2} \E{  (V(h) - V(0))^2 }, \quad h \in \mathbb{R}.
\]

In this contribution,  we shall discuss properties of
a functional closely related  to (\ref{wills}) and indirectly also related to the max-stable process $X$, namely
\[
	\mathcal{G}_V([S,T])=\E*{\inf_{t\in [S,T]}Z_V(t)},
	\quad 0\le S<T.
\]

In view of  {\it Brown-Resnick stationarity} and the fact that the infimum is a 1-homogeneous functional
(see the so-called {\it tilt-shift} formula, e.g., \cite{Htilt,Hrovje, hashorva2021shiftinvariant, hashorva2025cluster}),
it follows from \cite{Htilt} that
\bqn{\label{eq2C}
	\mathcal{G}_V([S,T]) = 	\mathcal{G}_V([0,T-S]).
}	
For notational simplicity, we will henceforth denote $\mathcal{G}_V([0,T])$ and $\mathcal{H}_V([0,T])$ by $\mathcal{G}_V(T)$ and $\mathcal{H}_V(T)$, respectively.
The constants  $\mathcal{G}_V(T)$ appear naturally in the tail asymptotics of inf-functionals of Gaussian
processes on short intervals,
\cite{DHL,DeK14,dkebicki2016extremes}.

For a general centered Gaussian process $V(t), t\inr$ with stationary increments,
a simple but useful explicit  upper bound is given by
\bqn{\label{eq:upperboundS}
\mathcal{G}_V(T) \le 2 \Psi(\sqrt{\gamma_V(T)}),\qquad T>0;}
see \Cref{PegAs} below. Clearly, \eqref{eq:upperboundS} implies that if
\bqn{\label{isra}
\limit{T} \gamma_V(T) = \IF,
}
then $\mathcal{G}_V(T)$ decays exponentially on the \(\gamma_V(T)\)-scale. This  is very different  from the linear growth of  $\mathcal{H}_V(T)$, when $	 \mathcal{H}_V>0$.\\

In view of \eqref{eq:upperboundS}, a natural problem is to determine whether under \eqref{isra}
\bqn{\label{eq:conL}
\lim_{T \to \IF} \frac{-\ln \mathcal{G}_V(T)}{\gamma_V(T)}=  A_V
\in (0, \infty) .
 }

It follows from  \Cref{PropKrzys} below that (\ref{eq:conL}) holds if
$V$ has  convex variogram $\gamma_V$.
In \Cref{thm:LDP} below we show \ks{the} above asymptotics for a large class of Gaussian processes with stationary increments $V$ including fBM $B_H$. For the latter, $A_{B_H}=1/2$ if $H\in [1/2,1)$ and $A_{B_H}$ is given by Gamma functions of $H$ if   $H\in (0,1/2)$, see \eqref{A_H1}.
This shows that  for  $H\in(0,1/2)$ the asymptotic behaviour of the expected infimum functional is
different from that when $H\in[1/2,1]$.

For a general Gaussian process $V(t), t\inr $ with variance function $\sigma_V^2$ we have
	\bqn{ \mathcal{G}_V(T) =\int_{\R} e^x \pk*{\inf_{t\in [0,T]}
		(\sqrt{2} V(t)- \sigma^2_V(t)) > x}  dx =: \int_{\R} e^x p_V(T,x) \, dx,
		\label{laj0}	
	}
which implies the following lower bound
	\bqn{\label{upB}
	 \mathcal{G}_V(T) \ge  e^{C}p_V(T,C), \quad C \in \R.
	}
The asymptotics of the persistence probability  $p_V(T,C)$ as $T \to \infty$ for general Gaussian processes has been studied in numerous contributions, see e.g., \cite{Molchan1999,aurzada2013persistence, AuS15,perB,perC,FFN21,aurzada2026,feldheim2025persistence} and the references therein.

In order to gain more intuition, we next discuss two tractable instances
investigated in \Cref{lemON} below, which illustrate how the covariance
structure influences both the decay of $\mathcal G_V(T)$ and its relation
to the corresponding persistence probabilities.
   \ks{Consider} first
$V(t)=B_1(t)=Zt$, $t\inr $ with $Z$ an $N(0,1)$ random variable (rv).
Writing
 $\Psi$ for the survival function of $Z$ and $\Phi=1-\Psi$ for its
distribution function, we have, for all $x<0$ and $T>0$, 
 \bqn{
	\label{nota0}
\mathcal{G}_{B_{1}}(T)= 2 \Psi(\sqrt{\gamma_{B_1}(T)}), \qquad
p_{B_1}(T,x)
		= \pk{\sqrt{2}Z T-T^2>x} =
		\Psi\left( \frac{x+ T^2}{\sqrt{2}T} \right),
} and hence
\bqn{\label{nota} p_{B_1}(T,x) \sim e^{-x/2} \mathcal{G}_{B_1}(T)/2, \qquad  T\to \IF.}
If $V=B_{1/2}$,  then
\bqn{\label{BMcase}
\mathcal{G}_{B_{1/2}}(T)
=
(2+T)\Psi\left(\sqrt{\frac{T}{2}}\right)
- \sqrt{\frac{T}{\pi}}\exp \left( -\frac{T}{4}\right)
}
implying
\bqn{ \label{boundA}
	\mathcal{G}_{B_{1/2}}(T) \sim \frac {\ks{4}} {\gamma_{B_{1/2}}(T)}  \Psi\left( \sqrt{\gamma_{B_{1/2}}(T)}\right)\sim  \frac{8e^{- T/4}}{T\sqrt{T \pi}}
		 , \quad T \to \IF,
}	
which shows that the decay of $\mathcal{G}_{B_{1/2}}(T) $ as $T\to \IF$ is very different from that of $\mathcal{G}_{B_{1}}(T)$.
The same dichotomy is also reflected in the corresponding persistence   probabilities, see below \Cref{ThmSlava}.  \\

Define below
$$\kTt :=\cov(V(t)\ks{-V(0)}, V(T)-V(t))= \gamma_V(T)- \gamma_V(t)- \gamma_V(T-t) .$$
For the first case $V(t)=tZ, t\ge 0$  we have
$$ \kTt  =t(T-t)\ge 0, \quad t\in [0,T], $$
while for the second case of Brownian motion
$$\kTt =0, \quad \forall t\in [0,T]. $$

If $\kTt$ remains non-negative for all $0 \le t\le T, T>0$, under assumption \eqref{isra},   \Cref{PropKrzys} below   implies the following sharp asymptotic bounds
\bqn{  \label{ausgez}
	\frac{4}{\gamma_V(T)}\Psi\left(\sqrt{\gamma_V(T)}\right)(1+o(1))
	 \le \mathcal{G}_{V}(T)    \le
	2\Psi\left(\sqrt{\gamma_V(T)}\right), \quad T \to \IF.
}

\kE{
In particular, \eqref{ausgez} holds for $V=B_H$ with $H\in [1/2,1]$ and when $H=1/2$ the lower bound is the exact asymptotics of $\mathcal{G}_{B_{1/2}}(T)$ as $T\to \IF$, while for $H=1$ the upper bound is the exact asymptotics of $\mathcal{G}_{B_{1}}(T)$ as $T\to \IF$.\\
We further investigate three classes of processes $V$ with bounded $\kTt$
for which the asymptotic behaviour of $\mathcal G_V(T)$ is similar to that
of Brownian motion, see \Cref{cor:BM-dominance-G}.\\
If $H\in (1/2,1)$, our result below shows that the upper bound in \eqref{ausgez} is asymptotically sharp, namely
\bqn{ \label{con1}
\mathcal{G}_{B_H}(T) \sim 2 \Psi( \sqrt{\gamma_{B_H}(T)}), \quad T\to \IF
}
and therefore a different asymptotic behaviour of $\mathcal{G}_{B_H}(T)$ is observed for $H\in (1/2,1)$ compared to that of $H=1/2$.
Moreover, for all $H\in (1/2,1]$ the asymptotics of the persistence exhibits a similar pattern. \\
For \(B_H\) with \(H>1/2\), the function \(\kappa_{B_H}(t,T)\) is non-negative
and grows in the interior of \([0,T]\), in a manner controlled by the variance
scale of \(B_H\), and suitable growth conditions on $\kTt$ away from the
endpoints allow us to establish \eqref{con1} for a large class of processes
satisfying \eqref{isra}, see \Cref{ThmSlava}.
}

If \(V\) is stationary with covariance function \(R_V\), then \eqref{isra}
cannot hold, since
\[
\gamma_V(t)=R_V(0)-R_V(t)
\]
is bounded. Moreover,
\[
\kTt
=
-R_V(0)-R_V(T)+R_V(t)+R_V(T-t),
\]
and therefore
\[
|\kTt|\le 4R_V(0),\qquad t\in[0,T],\ T>0.
\]
The associated process \(V_*(t)=V(t)-V(0)\), \(t\in\mathbb R\), has stationary
increments and the same variogram as \(V\), so $\kappa_{V_*}(t,T)$ is also
bounded and, in view of \Cref{PegAs} below, we have
\[
\mathcal{G}_V(T)=\mathcal{G}_{V_*}(T), \qquad T\ge0.
\]

For both $V$ and $V_*$ the variogram is bounded, hence the conjecture \eqref{con1} does not apply. In Theorem \ref{PropKE} we address a large class of ergodic stationary processes that satisfy Pickands-Berman conditions
(see (\ref{eqKrz}) and (\ref{eq:Berman})), for which we prove that
\bqn{ \label{as_G_st}
\mathcal{G}_V(T)
\sim
e^{\mathbb{E}\left\{\inf_{t\in[0,T]}\sqrt{2}V(t)-\sigma^2_V(0)\right\} }
\sim e^{-2\sigma_V(0) \sqrt{ \ln T}-\sigma_V^2(0)},
}
as $T\to\infty$.

\bigskip
{\bf The main contributions of the paper:}
\begin{enumerate}[(i)]
\item We pinpoint the crucial role played by the variogram $\gamma_V$ and the covariance $\kTt $ between $V(t)\ks{-V(0)}$ and the increment $V(T)-V(t)$ for $t \in [0,T]$.
Additionally, we obtain non-asymptotic bounds when $\kTt$ does not change sign.
\item We derive  in \ks{\Cref{thm:LDP}} logarithmic asymptotics  of   $\mathcal{G}_V(T)$  for  a large class of $V$ including the mixture of independent fBM's;
\item  For non-stationary $V$, under condition \eqref{isra} and an asymptotic growth condition on $\kTt$,
we refine the logarithmic asymptotic results by proving  \eqref{con1} for a large class of processes,
including fBM with Hurst parameter $H \in (1/2, 1]$. Further, we consider three classes of processes $V$ that exhibit similar behaviour to the Brownian motion case, where $\kTt$ is bounded;
\item For stationary \(V\), where condition \eqref{isra} cannot hold, we obtain
the exact asymptotics of \(\mathcal{G}_V(T)\) under the classical Pickands and
Berman conditions. Under the Berman condition alone, the corresponding
persistence probabilities decay faster than every power of $T$, on a
substantially different scale from \(\mathcal{G}_V(T)\).

\item Our results include new comparison principles for $\mathcal{G}_V(T)$ as well as results for the discrete analogues
of both $\mathcal{G}_V(T)$ and $p_V(T,C)$.
\end{enumerate}

{\bf Organisation}:\\
In  Section~2  we  discuss briefly the properties of   $\mathcal{G}_V(T)$ for general $V$ and then
establish in \Cref{PropKrzys} non-asymptotic bounds for the expected infimum.
Stationary \(V\) are discussed in \Cref{PropKE}; in this bounded-variogram
setting the large-\(\gamma_V(T)\) regime is absent, and ergodicity leads to a
different decay mechanism for \(\mathcal G_V(T)\).

Imposing further the classical Pickands and Berman conditions, we obtain the
exact asymptotics of \(\mathcal{G}_V(T)\). Under the Berman condition alone,
we also show that \(p_V(T,C)\) decays faster than every power of $T$.
In \Cref{ThmSlava}, for non-stationary $V$ satisfying  \eqref{isra}
and certain growth conditions on $\kTt$,
 we derive the exact asymptotics of both
$\mathcal{G}_V(T)$  and $p_V(T,C)$ as $T \to \infty$
for a broad class of log-normal \emph{Brown--Resnick stationary processes}, which includes in particular the fBM with Hurst parameter $H \in (1/2, 1]$.\\
Comparison principles for  $\mathcal{G}_V(T)$ are developed in \Cref{secC},
while \Cref{sec:dc} is devoted to the study of  discrete analogues of both $\mathcal{G}_V(T)$ and $p_V(T,C)$.   \\
All proofs, as well as explicit formulas for $\mathcal{G}_V(T)$   are deferred to \Cref{sec3}.

\section{Main Results}
As in the Introduction, we consider a separable centered Gaussian process $V(t)$, $t\in\R$,
with stationary increments, variance function $\sigma_V^2$, and variogram $\gamma_V$.

We first discuss   some general properties of the expected infimum functional.

\BEL  \label{PegAs}
If $V$ is separable, then for all $T>S $ we have:
\begin{enumerate}[(i)]
    \item \label{lemPeg:a} $\mathcal{G}_V([S,T]) = \mathcal{G}_V([0,T-S])$;
    \item \label{lemPeg:b} $\mathcal{G}_V(T):=\mathcal{G}_V([0,T]),T>0$ depends only on the variogram $\gamma_V$;
\item \label{lemPeg:c}   Setting $V_\star(t)=V(t)-V(0), t\ge 0$ we have
    \bqn{\label{taker}
	    \mathcal{G}_{V}(T)=\mathcal{G}_{V_\star}(T)\le
2 \Psi( \sqrt{\gamma_V(T)}), \quad T>0.
    }
\end{enumerate}
\EEL

We next derive sharp bounds for $\mathcal{G}_V(T)$ in the two cases where
\[
\kTt=\gamma_V(T)-\gamma_V(t)-\gamma_V(T-t),\qquad t\in[0,T],\ T>0,
\]
does not change sign.

\begin{proposition}\label{PropKrzys}  If  $V$ has continuous sample paths, then we have:
\begin{enumerate}[(i)]
\item
If $\kTt \ge 0$ for all $0\le t\le T$ and  $T>0$, then
\bqn{\label{lL1}
2\Psi\!\left(\sqrt{\gamma_V(T)}\right)
\ge
\mathcal{G}_{V}(T)
\ge
\mathcal{G}_{B_{1/2}}(2\gamma_V(T)).
}
Moreover, if $\limit{T}\gamma_V(T)=\infty$, then \eqref{ausgez} holds, and  \eqref{eq:conL} is satisfied with $A_V=1/2$.

\item
If $\kTt \le 0$ for all $0\le t\le T$ and all $T>0$, then
\bqn{\label{lL2}
\mathcal{G}_{V}(T)
\le
\mathcal{G}_{B_{1/2}}\!\left(2\max_{t\in[0,T]}\gamma_V(t)\right).
}
\end{enumerate}
\end{proposition}

\BRM \label{rem1}
\begin{enumerate}[(i)]
\item
If, in addition to the assumptions of \Cref{PropKrzys}, $V(0)=0$
almost surely, the argument therein yields the following bounds for every
$C_T<0$.
When $\kTt\ge0$ for all $0\le t\le T$ and $T>0$,
\[
p_{B_{1/2}}(2\gamma_V(T),C_T)
\le p_V(T,C_T)
\le p_{B_1}(\sqrt{2\gamma_V(T)},C_T).
\]
The upper bound follows by retaining only the endpoint $T$.
When $\kTt\le0$ for all $0\le t\le T$ and $T>0$,
\[
p_V(T,C_T)\le
p_{B_{1/2}}\!\left(2\max_{t\in[0,T]}\gamma_V(t),C_T\right).
\]
For an unanchored process these comparisons apply to $V_\star(t)=V(t)-V(0)$;
the identity $\mathcal G_V=\mathcal G_{V_\star}$ does not in general
extend to persistence probabilities.

\item If  $\gamma_V$ is convex, then $\kTt,t\in [0,T], T>0$ is non-negative. In this case, if further  \eqref{isra} holds, then \eqref{ausgez} confirms the logarithmic asymptotics in \eqref{eq:conL} with $A_V=1/2$. If $\gamma_V$ is concave, then $\kTt,t\in [0,T], T>0$ is non-positive.
\item If $\limit{T}\gamma_V(T)=L<\infty$, then the non-negativity of $\kTt, t\in [0,T],T>0$ forces the variogram $\gamma_V$ to be trivial, i.e., $\gamma_V(t)=0$ for all $t\ge0$. Hence a convex  non-trivial variogram $\gamma_V$ is necessarily unbounded. However,   there are \ks{also} genuine non-convex variograms $\gamma_V$ for which $\kTt,t\in (0,T), T>0$ is positive. For instance,  consider
\[
V(t)=B_H(t)+a\bigl(S(t)-S(0)\bigr),\qquad t\ge0, a>0,
\]
where $B_H, H\in (1/2,1]$ is an fBM independent of the Slepian process
\[
S(t)=B(t+1)-B(t),\qquad t\in\mathbb R
\]
with $B$ a standard Brownian motion.  Then,
\[
\sigma_V^2(t)=t^{2H}+2a^2\min(t,1),\qquad t\ge0,
\]
and for $a>0$ sufficiently small, $\kTt> 0$ for all $t\in (0,T)$ and all $T>0$.
\end{enumerate}
  \ERM

The upper bounds in \eqref{eq:upperboundS} and \eqref{lL2} do not force
$\mathcal{G}_V(T)$ to vanish as $T\to \IF$, if the variogram $\gamma_V(t),t\ge 0$ is bounded. This is in particular the case when  $V(t),t\inr$ is stationary. The special case   $V(t)=Z, t\inr$ with $Z$ a standard normal random variable yields
$$\gamma_V(T)=0, \qquad \mathcal{G}_V(T)=1, \qquad T>0.$$
 Another tractable strictly stationary but non-ergodic example is
$$V(t)=Z_1 \cos(t)+Z_2 \sin(t),\qquad t\inr , $$
with $Z_1,Z_2$ two independent $N(0,1)$ rv's and $T\ge  2\pi$. Again  we have that
$\mathcal{G}_V(T)$ is constant for all $T\ge 2\pi$ and $\gamma_V(t) =1- \cos(t),t\ge 0$ is bounded, see \Cref{lemON} below.

As shown below, if $V$ is ergodic, then  the expected infimum functional $\mathcal{G}_V(T)$
tends to zero as $T\to\infty$. Under the Pickands and Berman conditions
on $R_V(t):=\Cov(V(0),V(t))$, we also obtain its exact rate of decay.

\def\ve{\varepsilon}

\BT \label{PropKE}
If \(V(t), t\in\R\), is a centered stationary ergodic Gaussian process with
continuous sample paths and $
\Var(V(t))=\sigma^2\in(0,\infty)$,  then
\[
\lim_{T\to\infty}\mathcal G_V(T)=0.
\]
Additionally, if the Pickands condition
\bqn{\label{eqKrz}
R_V(0)- R_V(t) \ks{=}  A|t|^\alpha+ o(|t|^\alpha), \qquad t\to 0
}
is valid with \(\alpha\in(0,2]\), \(A>0\), and the Berman condition
\bqn{\label{eq:Berman}
\limit{T}\ln(T)R_V(T)= 0
}
holds, then \eqref{as_G_st} holds.
Under the Berman condition \eqref{eq:Berman}, for every $C\in\mathbb R$
and every $m>0$, we have
\begin{equation}\label{eq:persistenceA}
p_V(T,C)=o(T^{-m}),\qquad T\to\infty.
\end{equation}
\ET

\def\ve{\varepsilon}
\BRM
\begin{enumerate}[(i)]
\item The persistence bound \eqref{eq:persistenceA} uses only the Berman
condition \eqref{eq:Berman}; neither the local Pickands condition
\eqref{eqKrz} nor nonnegativity of the covariance is required.

\item
Following the subadditivity argument of \cite{LiS01}, if \(R_V(t)\ge 0\) for all
\(t\ge 0\), then for all \(C\in \R\) there exists
\begin{eqnarray}
b(C):=-\lim_{T\to\infty}  T^{-1}\ln \left(\pk*{\inf_{t\in[0,T]}(\sqrt2 V(t)-\sigma^2)>C}\right).\label{Shao}
\end{eqnarray}
The exact value of \(b(C)\) is known only for some particular processes; see
\cite{LiS01}. In the special case of \(C=-\sigma^2\), it is known that
\(b(C)>0\) is equivalent to \(\int_0^\infty R_V(t) dt<\infty\),
see \cite[Lem.~3.2]{perC}, while the
general case of \(C\in \R\) is studied in \cite{feldheim2025persistence}.
Persistence probabilities for processes with nonintegrable covariance
functions, under additional correlation or spectral assumptions, are studied
in \cite{DeS17,FFN21}.
This shows that persistence probabilities can exhibit a different asymptotic
behaviour compared to the expected infimum functional.
\end{enumerate}
\ERM

 We consider next the case of non-stationary $V$ satisfying further \eqref{isra}.
 
Our first result yields logarithmic asymptotics under regular variation
of the variance, and its proof relies on a Gaussian exponential tilt
determined by an explicit probability measure on $[0,1]$.

\begin{theorem}
\label{thm:LDP} Let $V(t),t\inr$ be a centered Gaussian process with stationary increments, continuous sample paths and
variance function $\sigma_V^2(t)$ and $\sigma_V(0)=0$.
Assume that $\sigma_V$ is regularly varying at infinity with index
$H\in(0,1)$.
Then
\[
\lim_{T\to\infty}
\frac{-\ln \mathcal G_V(T)}{\gamma _V(T)}
=
 A_{B_H},
\]
where
\begin{eqnarray}\label{A_H1}
A_{B_H}=\begin{cases}
        \frac{\Gamma\!\left(H+\frac12\right)\Gamma(1-2H)}
{\Gamma\!\left(\frac12-H\right)} \in ( 1/2,1),\quad &\text{ if } H \in(0,1/2),\\
        1/2,\quad &\text{ if } H \in [1/2,1\ks{)}.
    \end{cases}
\end{eqnarray}
\kE{Moreover, in the special case $V=B_H$ we have further
for all   $\aa{C}<0$
\bqn{\label{eq:conL2}
\lim_{T \to \IF} \frac{-\ln p_{B_H}(T,\aa{C})}{\gamma_{B_H}(T)}=  A_{B_H}.
 }
}
\end{theorem}

\BRM \begin{enumerate}[(i)]
	\item
Let $
V(t)=\sum_{i=1}^k a_i B_{H_i}(t), t\inr,$
where the \(B_{H_i}\)'s are independent fBM's with $H_i\in(0,1)$ and
$a_1,\ldots,a_k\in\mathbb R\setminus\{0\}$.
Then the assumptions
of \Cref{thm:LDP} are satisfied with
\[
\aa{H=\max\{H_1,\ldots, H_k\}.}
\]
Hence \Cref{thm:LDP} implies
\[
\lim_{T\to\infty}
\frac{-\ln \mathcal G_V(T)}{\gamma_V(T)}
=
A_{B_H}.
\]
\item \kE{In view of \eqref{upB}, under the assumptions of \Cref{thm:LDP}, we have, for
every \aa{\(C<0\)},
\[
\liminf_{T \to \IF}
\frac{-\ln p_V(T,C)}{\gamma_V(T)}
\ge A_{\aa{B_H}} .
\]
}
\end{enumerate}
\ERM

In order to obtain exact asymptotics of both \(\mathcal G_V\) and \(p_V\),
we impose conditions ensuring that the main contribution comes from the endpoint
\(T\). We require a growth bound on
$\kTt=\gamma_V(T)-\gamma_V(t)-\gamma_V(T-t)$ in the interior of the interval,
together with $\kappa_V(1,T)\to\infty$.
As shown in the proof, these conditions imply that the negative part of
$\kTt$ is asymptotically negligible; no separate sign condition is needed.
Under these hypotheses, the next result establishes \eqref{con1} for a large
class of non-stationary processes.

\BT
\label{ThmSlava}
Suppose that $V$ has continuous sample paths, $\sigma_V(0)=0$, and there exist constants $M > 1$ and $\epsilon > 0$ such that for all $T \geq 2M,  t \in [M, T/2]$
	\BQN
		\label{conv1}
		\kTt
		\geq {2} (1+\epsilon) \sigma_V(t) \sqrt{\ln t}.
	\EQN
	If further  
	\BQN	
		\label{conv2}
		\kappa_V(1,T)\longrightarrow\infty,\qquad T\to\infty,
	\EQN
then
	\bqn{ \label{con1:res}
\mathcal{G}_V(T)  \sim 2  \Psi( \sqrt{\gamma_V(T)}), \quad T\to \IF.
}
Moreover, for any threshold $C<0$
	 we have
	\BQN
		\label{eqpersistence}
		\quad \pk*{\inf_{t \in [0,T]} (\sqrt{2}V(t)-\sigma_V^2(t)) > C}
		\sim  \pk{\sqrt{2}V(T)-\sigma_V^2(T)>C}, \quad T \to \IF.
	\EQN
\ET
\BRM
\begin{enumerate}[(i)]
\item \kE{Both conditions
\eqref{conv1}  and \eqref{conv2} of \Cref{ThmSlava}  fail for   Brownian motion, since in this case $\kTt \equiv 0$ for all $t\in [0,T],T>0$.
In contrast, the conditions of \Cref{ThmSlava} are satisfied for
$
V(t)=\sum_{i=1}^k a_i B_{H_i}(t), t\inr,$
where the \(B_{H_i}\)'s are independent fBM's
\aa{with $\max\{H_1,\ldots, H_k\}>1/2$
and \(a_1,\ldots,a_k\in\mathbb \aa{R\setminus \{0\}}\)}}.
\item If $\sigma_V^2$ is convex,
 regularly varying at infinity with index $\alpha>1$ and $\sigma_V^2(0)=0$, then
 \eqref{conv1} and \eqref{conv2} are satisfied.
\item \kE{
The conditions of \Cref{ThmSlava} are stable under bounded uniformly continuous
perturbations of the variogram. More precisely, let
\[
V(t)=Y_1(t)+Y_2(t),\qquad t\in\mathbb R,
\]
where \(Y_1,Y_2\) are independent centered Gaussian processes with stationary
increments, continuous sample paths, and
\(\sigma_{Y_1}(0)=\sigma_{Y_2}(0)=0\). Suppose that \(Y_1\) satisfies
\eqref{conv1} and \eqref{conv2}. If the variogram \(\gamma_{Y_2}\) is bounded and
uniformly continuous on \([0,\infty)\), then \(V\) also satisfies
\eqref{conv1} and \eqref{conv2}. Hence the conclusions of \Cref{ThmSlava} hold for
\(V\).
In particular, this applies to
\[
V(t)=B_H(t)+Z(t)-Z(0),\qquad t\in\mathbb R,\quad H>1/2,
\]
where \(Z\) is a centered stationary Gaussian process with continuous sample
paths, independent of \(B_H\). Indeed, for
\(Y_2(t)=Z(t)-Z(0)\),
\[
\gamma_{Y_2}(t)=R_Z(0)-R_Z(t),
\]
which is bounded and uniformly continuous on \([0,\infty)\). No sign condition
on \(\kappa_{Y_2}\) is needed.
}
\item Let  \xx{$Z(s),s\ge 0$} be a centered stationary Gaussian process with continuous sample
 path\ks{s} and  covariance function
 $Cov(\xx{Z}(0),\xx{Z}(t))=R(t)\ge 0, \ t\ge 0$ and set
 $V(t)=\int_0^t \xx{Z}(s)ds, \ t\ge 0$.
 If $R$ is a regularly varying function at $\infty$ with  index $-\alpha \in(-1,0)$,
 then $\sigma^2_V(t)$ is a regularly varying function at $\infty$ with parameter $2-\alpha$ and convex, hence
 the conditions of Theorem \ref{ThmSlava} are satisfied.

 \item For an fBM $B_H$ with $H\in(0,1)$ and every fixed $C<0$,
\[
p_*(T,C):=\mathbb P\left\{\inf_{t\in[0,T]}B_H(t)>C\right\}
=T^{-(1-H)+o(1)},\qquad T\to\infty;
\]
see \cite{Molchan1999} and \cite[Thm.~12]{aurzada2018persistence}.
The sharper bound $p_*(T,C)\le cT^{-(1-H)}$, for some $c=c(H,C)>0$
and all large $T$, follows from the latter theorem for $H>1/2$, and
from the reflection principle for $H=1/2$.
Thus the polynomial decay of $p_*(T,C)$ differs from the exponential
decay of $p_{B_H}(T,C)$ on the scale $T^{2H}$.

\end{enumerate}
\ERM

\subsection{Comparison principles for $\mathcal{G}_{V}(T)$}
\label{secC}
We next establish a comparison between $\mathcal G_{V_1}(T)$ and
$\mathcal G_{V_2}(T)$, together with the corresponding relation for
$\mathcal H_{V_i}(T)$, $i=1,2$, extending slightly
\cite[Thm.~3.1]{debicki2017approximation}; these comparisons also yield
bounds for mixtures of independent fBM's and the interpolation inequalities
presented below. \\
 
\begin{theorem}\label{ThmKrzys}
Let   $V_i(t), t\inr,i=1,2$
be two centered  Gaussian processes with stationary increments, continuous sample paths,  and variogram $\gamma_{V_i},i=1,2$.
If further \kE{for some $  T>0$}
\bqn{ \eH{\gamma_{V_1}(t) \le \gamma_{V_2}(t), \quad \kE{\forall t\in [0,T]},
 }
 \label{gammas}
}
then we have
\bqn{\label{nachsitz}
\mathcal{G}_{V_1}(T)\ge \mathcal{G}_{V_2}(T), \quad
\mathcal{H}_{V_1}(T)\le \mathcal{H}_{V_2}(T).
}
\end{theorem}

\begin{corollary}
\label{cor:BM-dominance-G}
Let \( B(t),t\in\mathbb{R}\) be a standard Brownian motion independent of
\(Z(t),t\in\mathbb{R}\), which  is a centered stationary Gaussian process
with continuous sample paths,   \(\Var(Z(t))=\sigma^2\in(0,\infty)\) satisfying the  Pickands condition \eqref{eqKrz} for some \(A>0\), \(\alpha\in(0,2]\). Let further  \(V_c(t),t\in\mathbb{R}\) be a centered  Gaussian process with stationary
increments, continuous sample paths and variogram $\ks{\gamma_{V_c}(t)=} \sqrt{t^2+c\sin^2(t/2)}/2, c>0.$	
\begin{enumerate}[(i)]
\item
If \(Z\)
satisfies the Berman condition \eqref{eq:Berman}, then for $V(t)=B(t)+Z(t)-Z(0), t\in\mathbb{R}$
we have
\bqn{\label{eqH}
e^{-2\sigma\sqrt{\ln T}+O(1)}
\le
\frac{\mathcal G_V(T)}{\mathcal G_B(T)}
\le 1,
\qquad T\to\infty .
}
\item If $V=V_c$,  then \eqref{eqH} holds for  some \(\sigma\in(0,\infty)\)\ks{.}
\kE{
\item
Let $
V(t)=\int_0^t Z(s)\,ds, t\ge0,$  and suppose that $
R(t)=\Cov(Z(0),Z(t)) \ge0, t\ge0$. Assume further that
\[
a:=\int_0^\infty R(s)\,ds\in(0,\infty),
\qquad
m:=\int_0^\infty sR(s)\,ds<\infty
\]
and set $
g(t):=\int_t^\infty (s-t)R(s)\,ds, t\ge0$. If  \(t\mapsto g(|t|)\) is positive definite on
\(\mathbb R\) and  $
\lim_{t\to\infty}\ln(t)g(t)=0$,
then
we have
\[
1
\le
\frac{\mathcal G_V(T)}{\mathcal G_B(2aT)}
\le e^{2\sqrt m\,\sqrt{\ln T}+m+o(1)},
\qquad T\to\infty .
\]}
\end{enumerate}
\end{corollary}
\BRM
\label{rem:BM-critical-examples} \begin{enumerate}[(i)]
\item
Let $B(t_1,t_2,t_3), t_i\inr, i=1,2,3$ be the L\'evy Brownian motion. If $c>0$, then the process $V_{\aa{c}}(t):=B\left(t,\frac{\sqrt{c}}{2} \sin(t), \frac{\sqrt{c}}{2} \cos(t)\right), t\ge 0$,
has variogram $\sqrt{t^2+c\sin^2(t/2)}/2$, \xx{as supposed in Corollary \ref{cor:BM-dominance-G}}    .
\item \kE{
For the process $V$ in (iii) of \Cref{cor:BM-dominance-G} we have that $\kTt \in [0,m], t\in [0,T], T\ge 0$ and therefore it does not satisfy the assumptions of \Cref{ThmSlava}.
\item The assumptions in \(\Cref{cor:BM-dominance-G}\)(iii) are satisfied by
\[
Z_k(t)=B(t+k)-B(t),\qquad t\in\mathbb R, k>0.
\]
with
$a=
{k^2}/{2},
m=
{k^3}/{6}$ and
$
g(t)
=
{(k-t)_+^3}/{6}, t\ge0.
$
Note that  \(g(|t|)=k^3(1-|t|/k)_+^3/6,t\inr\)   is a scaled Askey covariance and
hence it is positive definite on \(\mathbb R\).
}
\end{enumerate}
\ERM

As a direct corollary of Theorem~\ref{ThmKrzys}, we obtain the following interpolation inequality for two fBM's.

\begin{corollary} \label{cor.comp}
Let $0<H_1< H< H_2\le 1$ and $T>0$. If $B_{H_1}$ and $B_{H_2}$ are two independent fBM's,  then setting
\[V(t)=\lambda^{1/2} T^{H-H_1}B_{H_1}(t)+(1-\lambda)^{1/2}T^{H-H_2}B_{H_2}(t),\qquad t\in[0,T],\]
where $\lambda=\frac{H_2-H}{H_2-H_1}$
we have
\[
\mathcal{G}_{B_H}(T)
\ge
\mathcal{G}_{V}(T)
\ge
\mathcal{G}_{B_{H_1}}\left(\lambda^{1/ (2H_1)}T^{H/H_1}\right)
\mathcal{G}_{B_{H_2}}\left((1-\lambda)^{1/ (2H_2)}T^{H/H_2}\right)
\]
and
\[
\mathcal{H}_{B_H}(T)
\le
\mathcal{H}_{V}(T)
\le
\mathcal{H}_{B_{H_1}}\left(\lambda^{1/ (2H_1)}T^{H/H_1}\right)
\mathcal{H}_{B_{H_2}}\left((1-\lambda)^{1/ (2H_2)}T^{H/H_2}\right).
\]
\end{corollary}

\begin{remark}\label{rem29}
In view of \eqref{eq:conL}, \Cref{cor.comp} yields that for every
$0<H_1<H<H_2\le 1$,
\[
A_{B_H}
\le
\lambda A_{B_{H_1}}+(1-\lambda)A_{B_{H_2}},
\qquad
\lambda=\frac{H_2-H}{H_2-H_1}.
\]
Thus, the map $H\mapsto A_{B_H}$ is convex on $(0,1]$. In particular, since $A_{B_{1/2}}=\frac12$, we obtain
\[
A_{B_H}
\le
\lambda A_{B_{H_1}}+\frac{1-\lambda}{2},
\qquad
0<H_1<H<\frac12,\qquad
\lambda=\frac{\frac12-H}{\frac12-H_1}.
\]
\end{remark}

 \subsection{Discrete constants }
 \label{sec:dc}

 Initially, the discrete Pickands constant
 $$ \mathcal{H}_V^\delta(T)=   \mathcal{H}_V([0,T]\cap  \delta \mathbb{Z}), \quad \delta>0
 $$
 \ks{was} defined in \cite{PickandsB} for $V=B_H, H\in (0,1]$.  Pickands' pioneering method of dealing with extremes of stationary Gaussian processes consists in a discretisation idea, where the discrete Pickands constants appear as the leading prefactors of the tail asymptotics of supremum of those processes. Pickands showed that the discrete constants converge to the continuous one $ \mathcal{H}_V^0$ when $\delta \to 0$ and therefore it was possible to \ks{analyse} the suprem\ks{a} of Gaussian processes on continuous time. That result has been extended for general Gaussian processes with stationary increments, see \cite{ZKE}. \\

If
\[
Z_V(t)\to0\quad\text{almost surely as }|t|\to\infty,
\]
then, for every $\delta\ge0$, the discrete Pickands constant satisfies
\[
\mathcal H_V^\delta
=\lim_{T\to\infty}\frac{\mathcal H_V^\delta(T)}{T}
\in(0,\infty),
\]
where $\delta=0$ denotes the continuous-time case; see, e.g.,
\cite{debicki2017approximation,ZKE}.

The expression for $\mathcal{H}_V^\delta $ for the particular case   $V= B_{1/2}$
has been calculated in \cite[Lem 5.16]{kabluchko2014limiting},  see also \cite{bisewski2025speed}. It is also of interest to study the rate of decrease of the discrete infimum constants
 $$\mathcal{G}_V^\delta(T) = \mathcal{G}_V([0,T]\cap  \delta \mathbb{Z})$$ as $T\to \IF$
   for a fixed constant $\delta>0$.

\BT
\label{ThmSlava2}
Under the assumptions of \Cref{ThmSlava}, for all \(\delta>0, C<0\),
\bqn{ \label{con12}
\mathcal{G}_V^\delta(T )
\sim
\mathcal{G}_V\left(\delta\left\lfloor T/\delta\right\rfloor\right)
\sim
2\Psi\left(
\sqrt{\gamma_V\left(\delta\left\lfloor T/\delta\right\rfloor\right)}
\right),
\quad T\to\infty
}
\bqn{ \label{eqpersistence2}
\pk*{\inf_{t \in [0, T] \cap \delta \mathbb{Z}}
(\sqrt{2}V(t)-\sigma_V^2(t)) > C}
\sim
p_V\left(\delta\left\lfloor T/\delta\right\rfloor,C\right),
\quad T \to \infty .
}
Further, for all  $\delta>0$
	\[
		\mathcal{G}_{B_{1/2}}([0,\delta N] \cap \delta \mathbb{Z}) \sim c_\delta\frac{e^{-\delta N/4} } {(\ks{\delta} N)^{3/2}}, \quad N \to \IF,
	\]
	where
\[
 c_\delta
 =\frac{2\delta}{\sqrt\pi}
 \exp\left\{2\sum_{k=1}^{\infty}
 \frac{e^{\delta k/4}}{k}
 \Psi\left(\sqrt{\delta k/2}\right)\right\}
 \longrightarrow\frac{8}{\sqrt\pi},\qquad\delta\downarrow0.
\]
Equivalently, we have 
\[
 c_\delta
 =\frac{2\delta}{\sqrt\pi}\sum_{k=0}^{\infty}
 e^{\delta k/4}
 \mathcal G_{B_{1/2}}([0,\delta k]\cap\delta\mathbb Z).
\]
	\ET

	\BRM
The comparisons in \Cref{ThmKrzys} remain valid when $[0,T]$ is replaced
by $[0,T]\cap\delta\mathbb Z$, whereas the discrete counterparts of the
interpolation inequalities in \Cref{cor.comp} require the lattice spacing
to be rescaled together with the time interval. More precisely,
self-similarity implies for all $a>0$, $H\in(0,1]$ and $\delta,T>0$
\[
\mathcal G_{aB_H}^{\delta}(T)
=\mathcal G_{B_H}^{a^{1/H}\delta}(a^{1/H}T),
\qquad
\mathcal H_{aB_H}^{\delta}(T)
=\mathcal H_{B_H}^{a^{1/H}\delta}(a^{1/H}T),
\]
so the same factor must multiply the lattice spacing and the horizon in
each term obtained from self-similarity.

 \ERM

\section{Further results and Proofs}
\label{sec3}
We first derive explicit formulas for $\mathcal G_V(T)$ in three tractable
cases, after which we proceed with the proofs of the results of Section~2.

\BEL \label{lemON}
For all \(T>0\),
\[
\mathcal G_{B_1}(T)=2\Psi\left(\sqrt{\gamma_{B_1}(T)}\right)
=2\Psi\left(\frac{T}{\sqrt2}\right).
\]
Moreover, \eqref{BMcase} holds for all \(T>0\), and \eqref{boundA} is
satisfied. Finally, if
\[
V(t)=Z_1\cos(t)+Z_2\sin(t),\qquad t\in\mathbb R,
\]
where \(Z_1,Z_2\) are independent \(N(0,1)\) rv's, then
\bqn{\label{eq:examp}
\mathcal{G}_{V}(T)=e^{-1}- 2\sqrt{\pi} \Psi(\sqrt{2}), \quad T\ge 2\pi.
}
\EEL

\prooflem{lemON} We begin by proving the formula for $\mathcal{G}_{B_1}(T)$.
Taking $W$ to be a standard normal rv and $c= \sqrt{2}\sigma_{B_1}(T)= \sqrt{2} T$, we obtain for all $T>0$
\bqny{
	\mathcal{G}_{B_{1}}(T)
	 &=& \E*{  e^{\min(0,  cW- c^2/2)}}\\
	 &=&\E*{ e^{  cW -c^2/2 }\ind{ cW - c^2/2 \le 0}}+ \pk{ W \ge c/2} \\
	 &=& 2\pk{ W\ge c/2}= 2 \Psi\Bigl( \frac{T}{\sqrt{2}}\Bigr).
}
Note in passing that the above argument shows that if $V(0)=0$, then
\bqn{ \label{boundBy}
    		\E*{\min\left(1, e^{
			\sqrt{2} V(t)-\sigma_V^2(t) }\right )} =
2\Psi\left(\sigma_V(t)/\sqrt{2}\right)   	, \quad t>0.
}

Next, we consider the case $H=1/2$.
Using that
\[
\mathcal{G}_{B_{1/2}}(T)
=
\E*{ e^{-\sqrt{2}\sup_{t\in[0,T]} (B(t)+t/\sqrt{2} ) }}
\]
and the explicit formula
\[
\frac{{\rm d}}{{\rm d}x} { \pk*{\sup_{t\in[0,T]} (B(t)+t/\sqrt{2} )\le x}}
=\sqrt{\frac{2}{\pi T}}\exp \left( - \frac{(x-T/\sqrt{2} )^2}{2T}\right)
 -\sqrt{2}e^{\sqrt{2}x}\Psi\left( \frac{x+T/\sqrt{2}}{\sqrt{T}} \right), \quad x\ge 0
\]
see, e.g., \cite{BD},  we obtain
\[
\mathcal{G}_{B_{1/2}}(T)=(2+T)\Psi\left(\sqrt{\frac{T}{2}}\right)
- \sqrt{\frac{T}{\pi}}\exp \left( -\frac{T}{4}\right),
\]
which confirms \eqref{BMcase}.

Letting  $\varphi=-\Psi'$ we have  the following asymptotic expansion for Mills ratio
\bqn{   \frac{\Psi(T)}{\varphi(T)}= \frac{1}{T}- \frac{1}{T^3}+
  \frac{3}{T^5} (1+ o(1)), \quad T\to \IF \label{Mills},
}
which follows from the inequalities
\bqn{
\sum_{k=0}^{2n-1}(-1)^k(2k-1)!!\,x^{-(2k+1)}
<
\frac{\Psi(x)}{\varphi(x)}
<
\sum_{k=0}^{2n}(-1)^k(2k-1)!!\,x^{-(2k+1)} .
\label{Mills2}
}
valid for all $x>0$ and all positive integers $n$.

 Consequently, we have

\bqny{
	\mathcal{G}_{B_{1/2}}(T) \sim
	\frac{8}{T}\Psi\left(\sqrt{\frac{T}{2}}\right)   	, \quad T \to \IF,
}	 	
establishing \eqref{boundA}. \\
Finally, we prove \eqref{eq:examp}. First, note that \(\sigma_V^2(t)=\Var(V(t))=1\) for all \(t\) and hence
\[
\mathcal{G}_V(T)
=
\E*{\inf_{t\in[0,T]} e^{\sqrt{2}V(t)-1}}
=
e^{-1}\E*{\exp\Bigl(\sqrt{2}\inf_{t\in[0,T]}V(t)\Bigr)}.
\]
Moreover $V(t)=R\cos(t-\Theta)$, where
$
R=\sqrt{Z_1^2+Z_2^2}$ with $\Theta$ uniformly distributed on $[0,2\pi]$ independent of $R$. Consequently, for all $T\ge 2\pi$ we have
\[
\inf_{t\in[0,T]}V(t)=-R
\]
implying thus
\[
\mathcal{G}_V(T)=e^{-1}\E*{e^{-\sqrt{2}R}}, \qquad T\ge 2\pi.
\]
Since the density of $R$ is given by
\[
f_R(r)=r e^{-r^2/2}, \qquad r>0,
\]
direct integration yields
\[
\E{e^{-\sqrt{2}R}}
=
1-2\sqrt{\pi}e\Psi(\sqrt{2}),
\]
which proves \eqref{eq:examp}.
\QED
\\

\prooflem{PegAs}
Proof of \Cref{lemPeg:a} and  \Cref{lemPeg:b}: When $V$ is stochastically continuous, the proof follows from the tilt-shift formula in \cite[Lem 4.5, Lem 7.1]{hashorva2021shiftinvariant}. For our case of a separable process $V$, we establish the assertions
directly as follows.
First, we note that  the stationarity of increments of $V$  yields
 $$
Cov(\sqrt{2}V(t), \sqrt{2}V(S)) =  Var(V(t)) + Var(V(S)) - Var(V(t-S)- V(0))$$
for all real $t,S$.
In the light of \cite[Lem.~6.1]{Htilt}, for all $t\inr$ and $T>S$ we obtain
\bqny{
   \mathcal{G}_{V}([S,T]) &=& \E*{\inf_{t \in [S,T]} e^{\sqrt{2}V(t)-\sigma^2_V(t)}}\\
   &=&  \E*{\inf_{t \in [S,T]}e^{\sqrt{2}V(S)-\sigma^2_V(S)}
   e^{ \sqrt{2}(V(t)- V(S))- (\sigma^2_V(t)-\sigma^2_V(S)) 	}}\\
   &=&  \E*{\inf_{t \in [S,T]}
  e^{  \sqrt{2}(V(t)- V(S))+Cov(\sqrt{2}(V(t)\ks{-V(S)}), \sqrt{2}V(S)) - (\sigma^2_V(t)-\sigma^2_V(S)) 	}}\\
   &=&   \E*{\inf_{t \in [0,T-S]} e^{ \sqrt{2}(V(t)- V(0))- 2\gamma_V(t)	} }\\
   &=&    \mathcal{G}_{V_\star}([0,T-S]),
}
where $V_\star(t)= V(t)-V(0)$.  The above (take now $S=0$)
combined with $\gamma_V(t)=\gamma_{V_\star}(t)$, yields
$$ \mathcal{G}_{V}(T):=\mathcal{G}_{V}([0,T]) = \mathcal{G}_{V_\star}([0,T]) .$$

Consequently,  $ \mathcal{G}_{V}(T)$ depends only on the variogram $\gamma_V$. \\

Proof of \Cref{lemPeg:c}: Since for all $T>0$
\begin{eqnarray*}
	\inf_{t\in[0,T]} (  \sqrt{2}V_\star(t)-2 \gamma_{V_\star}(t))
	&\le&
	\min(0, \sqrt{2}V_{\ks{\star}}(T)-\sigma^2_{V_{\ks{\star}}}(T))
\end{eqnarray*}
and
$$ \mathcal{G}_{V_\star}(T)\le \E*{
	e^{\min( 0, \sqrt{2}V_\star(T)-2 \gamma_{V_\star}(T)) 	}},
	$$
then the claim follows   from \eqref{boundBy}.
\QED
\\

\prooftheo{PropKE}
For each $n\in\mathbb N$, stationarity and ergodicity imply that, almost surely,
\[
\frac1T\int_0^T\ind{V(t)<-n}\,dt
\longrightarrow\mathbb P\{V(0)<-n\}>0.
\]
Consequently, $\inf_{t\in[0,T]}V(t)\to-\infty$ almost surely, and since
\[
0\le\inf_{t\in[0,T]}e^{\sqrt2V(t)-\sigma^2}
\le e^{\sqrt2V(0)-\sigma^2},
\qquad \mathbb E e^{\sqrt2V(0)-\sigma^2}=1,
\]
the dominated convergence theorem yields $\mathcal G_V(T)\to0$.

Suppose next that \eqref{eqKrz} and \eqref{eq:Berman} hold. In view of the
Berman condition, we have $r_V(t):=R_V(t)/\sigma^2<1$ for $t\ne0$, since
otherwise stationarity implies $V(kt)=V(0)$ almost surely for every positive
integer $k$, contradicting $R_V(kt)\to0$.
Let $M_T=\sup_{0\le t\le T}V(t)$, $\mu_T=\mathbb E M_T$, and let $d_T$
be a median of $M_T$. In the light of the Gaussian-maximum limit theorem
\cite[Thm.~12.3.5]{Led}, we obtain
\[
M_T-\sigma\sqrt{2\ln T}\longrightarrow0
\quad\text{in probability},
\]
and consequently $d_T=\sigma\sqrt{2\ln T}+o(1)$.
Using the one-sided Borell--TIS bounds, we further obtain
\[
|\mu_T-d_T|\le\sigma\sqrt{2\ln2},
\qquad
\mathbb P\{|M_T-\mu_T|>u\}
\le2e^{-u^2/(2\sigma^2)},\qquad u>0.
\]
It follows that, for every fixed $b>0$,
\[
\sup_{T>0}\mathbb E\exp\{b|M_T-d_T|\}<\infty.
\]
Since $M_T-d_T\to0$ in probability, uniform integrability (taking
$b>\sqrt2$ for the exponential below) yields
\[
\mu_T-d_T\to0,
\qquad
\mathbb E e^{-\sqrt2(M_T-d_T)}\to1.
\]
In view of the symmetry of the centered Gaussian process, we obtain
\[
\begin{split}
\mathcal G_V(T)
&=e^{-\sqrt2d_T-\sigma^2}
  \mathbb E e^{-\sqrt2(M_T-d_T)}\\
&\sim e^{-\sqrt2\mu_T-\sigma^2}
\sim e^{-2\sigma\sqrt{\ln T}-\sigma^2},
\end{split}
\]
which proves \eqref{as_G_st}.

Finally, in order to prove \eqref{eq:persistenceA} using only \eqref{eq:Berman},
put $a=-(\sigma^2+C)/(\sqrt2\sigma)$ and $q=\Phi(a+2)<1$, and then fix
$m>0$ and choose $K>m/(-\ln q)$. For sufficiently large $T$, let
\[
n_T=\lfloor K\ln T\rfloor,\qquad
\Delta_T=T/n_T,\qquad N_T=n_T+1,
\]
and
\[
\rho_T=\max\left(0,\max_{1\le k\le n_T}
                     \frac{R_V(k\Delta_T)}{\sigma^2}\right).
\]
Since $\ln\Delta_T\sim\ln T$, the Berman condition yields
$\rho_T=o(1/\ln T)$, while the event defining $p_V(T,C)$ implies
$-V(k\Delta_T)/\sigma<a$ for $0\le k\le n_T$.
In view of Slepian's inequality, its probability is bounded by the corresponding
probability for
\[
\eta_k=\sqrt{\rho_T}\,Z+\sqrt{1-\rho_T}\,Z_k,
\qquad 0\le k\le n_T,
\]
where $Z,Z_0,\ldots,Z_{n_T}$ are independent standard normal variables.
For $\rho_T>0$, we condition on $Z$ and split at
$-1/\sqrt{\rho_T}$, which yields, for all sufficiently large $T$
\[
\begin{split}
p_V(T,C)
&\le\mathbb P\{\eta_k\le a,\ 0\le k\le n_T\}\\
&\le\Phi\left(\frac{a+1}{\sqrt{1-\rho_T}}\right)^{N_T}
     +\mathbb P\{Z<-1/\sqrt{\rho_T}\}\\
&\le q^{N_T}+e^{-1/(2\rho_T)}
=o(T^{-m}).
\end{split}
\]
If $\rho_T=0$, the comparison vector is independent and the same conclusion
follows from $p_V(T,C)\le\Phi(a)^{N_T}\le q^{N_T}$.
\QED

\proofprop{PropKrzys} In view of \Cref{PegAs} we can assume without loss of generality that $\sigma_V(0)=0$ implying that $\gamma_V(t)=\sigma_V^2(t)/2$ for all $t\ge 0$. Next,  for all $0\le s\le t<\IF$ \ks{such that $\sigma_V^2(s)\leq \sigma_V^2(t)$}
\[
\Cov(V(s),V(t))-\Cov(B_{1/2}(\sigma_V^2(s)),B_{1/2}(\sigma_V^2(t)))
=\kappa_V(s,t) .
\]

\medskip
(i) If $\kappa_V(s,t)\ge 0$, then \ks{$\sigma_V^2(t), t\ge 0$ is nondecreasing. Hence} $$\Cov(V(s),V(t))\ge \Cov(B_{1/2}(\sigma_V^2(s)),B_{1/2}(\sigma_V^2(t))).$$   Moreover,  $\{\sigma_V^2(t):t\in[0,T]\}=[0,\sigma_V^2(T)]$. In view of Slepian's inequality, for $x\le0$ we obtain
\[
\pk*{\inf_{t\in[0,T]}(\sqrt2V(t)-\sigma_V^2(t))>x}
\ge
\pk*{\inf_{u\in[0,\sigma_V^2(T)]}(\sqrt2B_{1/2}(u)-u)>x}.
\]
Integration yields $\mathcal{G}_{V}(T)\ge \mathcal{G}_{B_{1/2}}(2\gamma_V(T))$.
Using further   \eqref{taker} the claim follows easily.

\medskip
(ii) Assume that \(\kappa_V(s,t)\le0\) for all \(0\le s\le t\).
Put
\[
\theta(t):=\sigma_V^2(t),\qquad
m_T:=\max_{t\in[0,T]}\theta(t).
\]
If \(m_T=0\), then \(V(t)=0\) a.s. on \([0,T]\), and the claim is trivial.
Assume therefore that \(m_T>0\).

Since $\theta$ is continuous and $\theta(0)=0$, its range on $[0,T]$
is $[0,m_T]$. For \(u\in[0,m_T]\), define the
first hitting time
\[
\tau(u):=\inf\{t\in[0,T]:\theta(t)=u\}.
\]
Then \(\tau(u)\) is well defined, \(\theta(\tau(u))=u\), and \(\tau\) is
nondecreasing. Indeed, if \(0\le u_1<u_2\le m_T\), then before hitting level
\(u_2\), the continuous function \(\theta\) must hit level \(u_1\).

Let
\[
0=u_0<u_1<\cdots<u_n=m_T
\]
be arbitrary. Set \(t_i=\tau(u_i)\). Then \(t_i\le t_j\) whenever \(i\le j\),
and \(\theta(t_i)=u_i\). For \(i\le j\), using \(\kappa_V(t_i,t_j)\le0\), we get
\[
\begin{aligned}
&\Cov\bigl(\sqrt2 V(t_i),\sqrt2 V(t_j)\bigr)
-
\Cov\bigl(\sqrt2 B_{1/2}(u_i),\sqrt2 B_{1/2}(u_j)\bigr)
\\
&\qquad
=
2\Cov(V(t_i),V(t_j))-2u_i
\\
&\qquad
=
2\kappa_V(t_i,t_j)
\le0.
\end{aligned}
\]
Moreover, we have 
\[
\Var(\sqrt2 V(t_i))=2u_i
=
\Var(\sqrt2 B_{1/2}(u_i)).
\]
Consequently, an application of Slepian's inequality yields, for every
$x\in\mathbb R$ 
\[
\mathbb P\left\{
\inf_{t\in[0,T]}
\bigl(\sqrt2V(t)-\sigma_V^2(t)\bigr)>x
\right\}
\le
\mathbb P\left\{
\min_{0\le i\le n}
\bigl(\sqrt2B_{1/2}(u_i)-u_i\bigr)>x
\right\}.
\]

Now take a sequence of partitions whose mesh tends to zero. Since the process
\(u\mapsto \sqrt2B_{1/2}(u)-u\) has continuous sample paths and its
infimum over $[0,m_T]$ has a continuous distribution, we obtain
\[
\mathbb P\left\{
\inf_{t\in[0,T]}
\bigl(\sqrt2V(t)-\sigma_V^2(t)\bigr)>x
\right\}
\le
\mathbb P\left\{
\inf_{u\in[0,m_T]}
\bigl(\sqrt2B_{1/2}(u)-u\bigr)> x
\right\}.
\]
Since
\[
m_T=\max_{t\in[0,T]}\sigma_V^2(t)
=
2\max_{t\in[0,T]}\gamma_V(t),
\]
we conclude that
\[
\mathcal G_V(T)
\le
\mathcal G_{B_{1/2}}
\left(
2\max_{t\in[0,T]}\gamma_V(t)
\right),
\]
which proves \eqref{lL2}.
\QED
\\

\prooftheo{thm:LDP}
Writing $\beta=2H$ and $s_T=\sigma_V^2(T)$, we first record two
consequences of continuity, stationary increments and regular variation:
\begin{equation}\label{eq:tilt-rv-uniform}
 \sup_{0\le u\le1}\left|
 r_T(u)-u^\beta\right|\longrightarrow0,
 \qquad r_T(u):=\frac{\sigma_V^2(Tu)}{s_T},
\end{equation}
and
\begin{equation}\label{eq:tilt-mean-negligible}
 m_T:=\mathbb E\sup_{0\le t\le T}V(t)=o(s_T).
\end{equation}
For the first assertion, the uniform convergence theorem applies on
$[\delta,1]$ for each $\delta>0$, while on $[0,\delta]$ Potter's bound yields
for any $\kappa\in(0,\beta)$
\[
 \limsup_{T\to\infty}\sup_{0\le u\le\delta}r_T(u)
 \le C\delta^{\beta-\kappa}.
\]
Here the portion $Tu\le T_0$ is bounded by
$\sup_{0\le t\le T_0}\sigma_V^2(t)/s_T\to0$, so letting $\delta\downarrow0$
proves \eqref{eq:tilt-rv-uniform}.

To prove \eqref{eq:tilt-mean-negligible}, put $n=\lceil T\rceil$ and
\[
 M_k=\sup_{0\le u\le1}|V(k+u)-V(k)|,\qquad 0\le k<n.
\]
In view of stationary increments, these random variables have the same
distribution, while continuity and Fernique's theorem imply
$\mathbb E M_0<\infty$ and Borell's inequality yields
\[
 \mathbb P\{M_k>\mathbb E M_0+x\}
 \le\exp\{-x^2/(2v_1)\},\qquad x>0,
 \qquad v_1:=\sup_{0\le u\le1}\sigma_V^2(u)>0.
\]
A union bound and integration therefore yield
$\mathbb E\max_{k<n}M_k=O(\sqrt{\ln(T+2)})$,
and the standard bound for a finite Gaussian maximum implies
\[
 \mathbb E\max_{0\le k\le n}|V(k)|
 \le \sqrt{2\ln(2n+2)\max_{0\le k\le n}\sigma_V^2(k)}
 =O\!\left(\sigma_V(T)\sqrt{\ln(T+2)}\right),
\]
where \eqref{eq:tilt-rv-uniform} and regular variation imply the last
bound. Since
\[
 \sup_{0\le t\le T}|V(t)|
 \le \max_{0\le k\le n}|V(k)|+\max_{0\le k<n}M_k,
\]
we obtain
$m_T=O(\sigma_V(T)\sqrt{\ln(T+2)})=o(s_T)$, because $H>0$.

We next choose a probability measure $q_H$ on $[0,1]$ whose potential
\[
 U_H(t):=\int_0^1|t-s|^{2H}q_H(ds)
\]
satisfies
\begin{equation}\label{eq:tilt-potential}
 \sup_{0\le t\le1}U_H(t)
 =\int_0^1U_H(t)q_H(dt)
 =\int_0^1t^{2H}q_H(dt)=:c_H.
\end{equation}
For $H\ge1/2$, take $q_H=(\delta_0+\delta_1)/2$, so that $c_H=1/2$ and
\[
 U_H(t)=\frac{t^{2H}+(1-t)^{2H}}2\le\frac12,
\]
with equality at both endpoints.

For $H<1/2$, put $a=1/2-H$ and take $q_H$ to have the
Beta$(a,a)$ density, i.e.,
\[
 q_H(ds)=\frac{s^{a-1}(1-s)^{a-1}}{\mathrm B(a,a)}\,ds,
 \qquad 0<s<1,
\]
where $\mathrm B$ denotes Euler's beta function. Its potential is constant;
indeed, for $0<t<1$, the substitution
$s=tv/(1-t+tv)$ yields
\[
 \int_0^t(t-s)^{2H-1}s^{a-1}(1-s)^{a-1}\,ds
 =\mathrm B(a,2H)[t(1-t)]^{-a}.
\]
The corresponding integral over $(t,1)$ has the same value by symmetry, and
furthermore
\[
 \int_0^1\int_0^1 2H|t-s|^{2H-1}\,dt\,q_H(ds)
 =\int_0^1\bigl(s^{2H}+(1-s)^{2H}\bigr)q_H(ds)<\infty.
\]
Using Fubini's theorem, we conclude that $U_H$ is absolutely continuous on
$[0,1]$, and its derivative is zero almost everywhere in view of the equality of
the derivative contributions from $(0,t)$ and $(t,1)$; hence
$U_H$ is constant, and consequently
\begin{equation}\label{eq:beta-id}
 U_H(t)=\int_0^1s^{2H}q_H(ds)
 =\frac{\Gamma(H+\frac12)\Gamma(1-2H)}
        {\Gamma(\frac12-H)}=:c_H,
 \qquad 0\le t\le1.
\end{equation}
Furthermore, $1/2<c_H<1$, since $s<s^{2H}<1$ on $(0,1)$ and
$q_H$ is symmetric about $1/2$. This proves
\eqref{eq:tilt-potential} in both cases, with $c_H=A_{B_H}$ as specified
in \eqref{A_H1}.

Fix this measure $q_H$ and define
\[
 X_T(t)=V(Tt),\qquad
 Z_T=\int_0^1X_T(s)q_H(ds),\qquad
 c_T=\int_0^1r_T(s)q_H(ds),
\]
\[
 U_T(t)=\int_0^1r_T(|t-s|)q_H(ds),\qquad
 D_T=\int_0^1U_T(t)q_H(dt).
\]
In view of stationary increments and $V(0)=0$, we obtain
\[
 b_T:=\operatorname{Var}(Z_T)=s_T(c_T-D_T/2),
 \qquad
 2\operatorname{Cov}(X_T(t),Z_T)-\operatorname{Var}(X_T(t))
 =s_T(c_T-U_T(t)).
\]
Writing the original integrand as
\[
 e^{\sqrt2Z_T}\exp\left\{\inf_{0\le t\le1}
 \bigl[\sqrt2(X_T(t)-Z_T)-s_Tr_T(t)\bigr]\right\}
\]
and applying Gaussian exponential tilting yields the exact identity
\begin{equation}\label{eq:tilt-infimum-identity}
 \mathcal G_V(T)
 =e^{-b_T}\mathbb E\exp\left\{
 \inf_{0\le t\le1}
 \left[\sqrt2(X_T(t)-Z_T)+s_T(c_T-U_T(t))\right]
 \right\}.
\end{equation}
Using Jensen's inequality on the right-hand side, together with
$\mathbb E Z_T=0$ and Gaussian symmetry, we obtain
\[
 \ln\mathcal G_V(T)
 \ge s_T\left(\frac{D_T}{2}-\sup_{0\le t\le1}U_T(t)\right)
      -\sqrt2m_T.
\]
For the opposite bound, averaging the original infimum against $q_H$
yields directly
\[
 \mathcal G_V(T)
 \le\mathbb E\exp\{\sqrt2Z_T-s_Tc_T\}
 =\exp\{-s_TD_T/2\}.
\]
In the light of \eqref{eq:tilt-rv-uniform}, $U_T\to U_H$ uniformly, so
\eqref{eq:tilt-potential} implies $D_T\to c_H$ and
$\sup U_T\to c_H$, and hence, together with \eqref{eq:tilt-mean-negligible}, the
last two bounds yield
\[
 \lim_{T\to\infty}\frac{\ln\mathcal G_V(T)}{s_T}
 =-\frac{c_H}{2}.
\]
Since $\gamma_V(T)=s_T/2$, this proves the first assertion.

To prove the persistence assertion for $V=B_H$, let $C<0$ and
$\varepsilon=T^{-H}$, so that self-similarity implies
\[
 p_{B_H}(T,C)=\mathbb P\left\{
 \sqrt2\varepsilon B_H(u)-u^{2H}>C\varepsilon^2,
 \quad 0\le u\le1\right\}.
\]
With $Z=\int_0^1B_H(s)q_H(ds)$, \eqref{eq:tilt-potential} implies
\[
 \operatorname{Var}(Z)=c_H/2,
 \qquad
 2\operatorname{Cov}(B_H(u),Z)=u^{2H}+c_H-U_H(u)\ge u^{2H}.
\]
Define the tilted probability measure $\mathbb Q_\varepsilon$ by
\[
 L_\varepsilon:=\frac{d\mathbb Q_\varepsilon}{d\mathbb P}
 =\exp\left\{\frac{\sqrt2Z}{\varepsilon}
              -\frac{c_H}{2\varepsilon^2}\right\}.
\]
Under $\mathbb Q_\varepsilon$, the process $B_H$ has its original
covariance and mean
$\sqrt2\operatorname{Cov}(B_H(\cdot),Z)/\varepsilon$.
Writing $E_\varepsilon$ for the event defining $p_{B_H}(T,C)$ and
\[
 A_\varepsilon=\left\{
 \inf_{0\le u\le1}B_H(u)>C\varepsilon/\sqrt2\right\},
\]
we therefore have
$\mathbb Q_\varepsilon(E_\varepsilon)\ge\mathbb P(A_\varepsilon)$.
For every $\delta>0$, H\"older's inequality yields
\[
 \mathbb Q_\varepsilon(E_\varepsilon)
 \le\mathbb P(E_\varepsilon)^{\delta/(1+\delta)}
       \bigl(\mathbb E L_\varepsilon^{1+\delta}\bigr)^{1/(1+\delta)}
 =\mathbb P(E_\varepsilon)^{\delta/(1+\delta)}
       \exp\left\{\frac{\delta c_H}{2\varepsilon^2}\right\}.
\]
Consequently, we obtain 
\[
 p_{B_H}(T,C)\ge\mathbb P(A_\varepsilon)^{1+1/\delta}
 \exp\left\{-\frac{(1+\delta)c_H}{2\varepsilon^2}\right\}.
\]
To bound $\mathbb P(A_\varepsilon)$ and set   $d_C=-C/(2\sqrt2)>0$;
using Gaussian symmetry and self-similarity, we obtain
\[
 \mathbb P(A_\varepsilon)
 \ge\mathbb P\left\{\sup_{0\le u\le1}B_H(u)\le d_C\varepsilon\right\}
 =\mathbb P\left\{\sup_{0\le u\le\varepsilon^{-1/H}}B_H(u)\le d_C\right\}.
\]
The factor $1/2$ in the definition of $d_C$ ensures the strict inequality
in $A_\varepsilon$. In accordance with \cite[Thm.~12]{aurzada2018persistence}, whose
lower bound holds for every $H\in(0,1)$, we obtain
\[
 \mathbb P(A_\varepsilon)
 \ge c\,\varepsilon^{(1-H)/H}
          |\ln\varepsilon|^{-1/(2H)}
\]
for sufficiently small $\varepsilon>0$, with $c>0$ depending on $H$ and
$C$. It follows that $\varepsilon^2\ln\mathbb P(A_\varepsilon)\to0$, and consequently
\[
 \liminf_{\varepsilon\downarrow0}\varepsilon^2
 \ln p_{B_H}(T,C)\ge-\frac{(1+\delta)c_H}{2}.
\]
Letting $\delta\downarrow0$ and recalling
$\varepsilon^{-2}=T^{2H}=2\gamma_{B_H}(T)$ yields
\[
 \limsup_{T\to\infty}
 \frac{-\ln p_{B_H}(T,C)}{\gamma_{B_H}(T)}\le c_H.
\]
The reverse inequality follows from
$p_{B_H}(T,C)\le e^{-C}\mathcal G_{B_H}(T)$ and the first assertion,
which completes the proof.
\QED

\prooftheo{ThmSlava}
In view of \eqref{conv2}, we have $\kappa_V(1,T)\to\infty$ and since
\[
\kappa_V(1,T)\le\sigma_V(1)\sigma_V(T)-\sigma_V^2(1),
\]
it follows that $\sigma_V(T)\to\infty$. Moreover, $\sigma_V(t)>0$
for every $t>0$, since otherwise stationary increments would imply
$\sigma_V(kt)=0$ for all positive integers $k$, contradicting this limit.

We first derive two consequences of the hypotheses. Write
$\gamma=\gamma_V$ and $\sigma=\sigma_V$ and define for $r\in\mathbb R$
\[
f(r)=\int_{-1/2}^{1/2}\gamma(r+u)\,du,
\qquad b(r)=\gamma(r)-f(r).
\]
By continuity and \eqref{conv2}
\[
f'(r)=\gamma(r+1/2)-\gamma(r-1/2)
=\kappa_V(1,r+1/2)+\gamma(1)\longrightarrow\infty,
\qquad r\to\infty.
\]
For $0\le u\le1/2$ stationary increments implies 
\[
D_u(r):=\gamma(r+u)+\gamma(r-u)-2\gamma(r)
=\operatorname{Cov}(V(r+u)-V(r),V(u)).
\]
Consequently, we have 
\[
b(r)=-\int_0^{1/2}D_u(r)\,du,
\qquad
\|b\|_\infty\le\int_0^{1/2}\sigma^2(u)\,du<\infty.
\]
Moreover, the difference of two length-$u$ increments can be regrouped as
the difference of two length-$|h|$ increments. Cauchy--Schwarz therefore yields
\[
|D_u(r+h)-D_u(r)|\le2\sigma(u)\sigma(|h|),
\]
and hence
\[
|b(r+h)-b(r)|\le K\sigma(|h|),
\qquad K=2\int_0^{1/2}\sigma(u)\,du.
\]
Thus $b$ is bounded and uniformly continuous. For every fixed
$0<\delta<L<\infty$, since $f'$ is eventually positive, we obtain
\begin{equation}\label{conv2:uniform}
\begin{split}
\inf_{t\in[\delta,L]}\kappa_V(t,T)
&\ge\delta\inf_{r\in[T-L,T]}f'(r)
-2\|b\|_\infty-\sup_{t\in[\delta,L]}\gamma(t)\\
&\longrightarrow\infty,\qquad T\to\infty.
\end{split}
\end{equation}
In particular, \eqref{conv2} is equivalent to divergence of
$\kappa_V(t,T)$ uniformly on every compact subinterval of $(0,\infty)$.

For each fixed $\tau\in(0,M)$, eventual monotonicity of $f$ also implies,
for all sufficiently large $T$ and $0\le t\le\tau$
\[
\kappa_V(t,T)
=f(T)-f(T-t)+b(T)-b(T-t)-\gamma(t)
\ge-K\sigma(t)-\gamma(t).
\]
On $[\tau,M]$, \eqref{conv2:uniform} gives eventual nonnegativity, while
\eqref{conv1} does so on $[M,T/2]$. Using the symmetry
$\kappa_V(t,T)=\kappa_V(T-t,T)$ for the other half of the interval, we get
\[
\limsup_{T\to\infty}\sup_{t\in[0,T]}(-\kappa_V(t,T))_+
\le\sup_{t\in[0,\tau]}(K\sigma(t)+\gamma(t)).
\]
Letting $\tau\downarrow0$ proves
\begin{equation}\label{conv3}
\lim_{T\to\infty}\sup_{t\in[0,T]}(-\kappa_V(t,T))_+=0.
\end{equation}

Put
\[
S_T=\sigma_V(T),\qquad q_T(t)=\frac{R_V(t,T)}{S_T^2},\qquad
X_T(t)=\sqrt2\bigl(V(t)-q_T(t)V(T)\bigr),
\]
\[
W_T(t)=X_T(t)+\kappa_V(t,T),\qquad
M_T(x)=\inf_{t\in[0,T]}\bigl(W_T(t)+\sqrt2xq_T(t)\bigr),
\]
where $R_V$ is the covariance function of $V$. Let next 
\[
b_T=\frac{\varphi(S_T/\sqrt2)}{S_T},\qquad
h_T(x)=e^{-x/\sqrt2-x^2/(2S_T^2)}.
\]
The centered Gaussian bridge $X_T$ is independent of $V(T)$, and therefore
conditioning on $V(T)=S_T^2/\sqrt2+x$ yields
\[
\begin{aligned}
\frac{\mathcal G_V(T)}{b_T}
 &=\int_{\mathbb R}\mathbb E[e^{M_T(x)}]h_T(x)\,dx,\\
\frac{p_V(T,C)}{b_T}
 &=\int_{\mathbb R}\mathbb P\{M_T(x)>C\}h_T(x)\,dx,
\end{aligned}
\]
We first prove
\begin{equation}\label{probconv2}
L_T:=\inf_{t\in[0,T]}W_T(t)\stackrel{P}{\longrightarrow}0.
\end{equation}
Gaussian regression and stationary increments imply
\[
\operatorname{Var}(X_T(t)-X_T(s))
=2\sigma_V^2(t-s)
 -2\frac{(R_V(t,T)-R_V(s,T))^2}{S_T^2}
\le2\sigma_V^2(t-s).
\]
Since $X_T(0)=X_T(T)=0$, we also have
\[
\operatorname{Var}(X_T(t))
\le2\min\{\sigma_V^2(t),\sigma_V^2(T-t)\},\qquad 0\le t\le T.
\]
Write
\[
m_h=\mathbb E\sup_{0\le t\le h}\sqrt2V(t)<\infty.
\]
Continuity and integrability of compact Gaussian suprema imply
$m_h\to0$ as $h\downarrow0$, while for every interval $I=[u,u+h]\subset[0,T]$,
the Sudakov--Fernique inequality and stationary increments yield
$\mathbb E\sup_{t\in I}X_T(t)\le m_h$.
Consequently, using Gaussian symmetry and the Borell--TIS inequality,
we obtain  for $h>0$
\[
\mathbb P\left\{\inf_{t\in I}W_T(t)<0\right\}
\le
\exp\left\{
-\frac{\bigl(\inf_{t\in I}\kappa_V(t,T)-m_h\bigr)_+^2}
 {4\sup_{t\in I}\min\{\sigma_V^2(t),\sigma_V^2(T-t)\}}
\right\}.
\]
This bound holds directly also when the bridge is deterministic.

In order to use the same estimates at both endpoints, we first note that
stationary increments imply
\[
q_T(T-t)=1-q_T(t),\qquad
\kappa_V(T-t,T)=\kappa_V(t,T).
\]
The process $V(T)-V(T-t)$, $0\le t\le T$, has the same law as $V(t)$,
and since its bridge is $-X_T(T-t)$, Gaussian symmetry yields
\[
\{W_T(T-t):0\le t\le T\}
\stackrel{d}{=}\{W_T(t):0\le t\le T\}.
\]

To bound the middle of the interval, we observe that stationary increments,
together with the triangle inequality in $L^2$ imply
\[
|\sigma_V(u)-\sigma_V(v)|
\le\sigma_V(u-v)\le\sup_{0\le r\le1}\sigma_V(r)
\qquad\text{if }|u-v|\le1.
\]
In the light of $\sigma_V(k)\to\infty$, this inequality implies
\[
\frac{\inf_{t\in[k,k+1]}\sigma_V(t)}
     {\sup_{t\in[k,k+1]}\sigma_V(t)}\longrightarrow1.
\]
For integers $k\ge M$ with $k<T/2$ let
\[
J_{k,T}=[k,\min(k+1,T/2)],\qquad
v_{k,T}=\sup_{t\in J_{k,T}}\sigma_V(t).
\]
Using \eqref{conv1} together with the preceding ratio, we obtain uniformly in $T$
\[
\inf_{t\in J_{k,T}}\kappa_V(t,T)
\ge2(1+\epsilon)\inf_{t\in J_{k,T}}\sigma_V(t)\sqrt{\ln k}
\ge2(1+\epsilon/2)v_{k,T}\sqrt{\ln k}
\]
for all sufficiently large $k$. Since $v_{k,T}\ge\sigma_V(k)\to\infty$,
it follows, after increasing the lower bound on $k$ that
\[
\inf_{t\in J_{k,T}}\kappa_V(t,T)-m_1
\ge2(1+\epsilon/4)v_{k,T}\sqrt{\ln k}.
\]
Since each block has length at most one, its mean bound is at most $m_1$
and its variance denominator is at most $4v_{k,T}^2$ and hence applying
the interval bound above with
$\eta=(1+\epsilon/4)^2-1>0$ yields
\[
\mathbb P\left\{\inf_{t\in J_{k,T}}W_T(t)<0\right\}
\le k^{-1-\eta}.
\]
The last block is also covered by this estimate, even when its length is
less than one. Consequently, for every sufficiently large integer $L$ and all $T>2L$,
summing the blocks with $k\ge L$, $k<T/2$, and using reflection symmetry yields
\[
\mathbb P\left\{\inf_{t\in[L,T-L]}W_T(t)<0\right\}
\le2\sum_{k=L}^{\infty}k^{-1-\eta}.
\]

For fixed $0<\tau<L$, \eqref{conv2:uniform} implies
$\inf_{t\in[\tau,L]}\kappa_V(t,T)\to\infty$, while the mean bound $m_L$
and the variance bound $2\sup_{0\le t\le L}\sigma_V^2(t)$ are independent of $T$.
In view of these bounds, the same interval estimate yields
\[
\mathbb P\left\{\inf_{t\in[\tau,L]}W_T(t)<0\right\}\to0,
\]
and reflection yields the analogous conclusion on $[T-L,T-\tau]$.
Finally, set
\[
a_T=\sup_{0\le t\le T}(-\kappa_V(t,T))_+\to0,
\]
where the convergence follows from \eqref{conv3}.
For each $d>0$ and all sufficiently large $T$ such that $a_T<d/2$,
applying Markov's inequality together with the bound on the expected
supremum yields
\[
\begin{aligned}
\mathbb P\left\{\inf_{t\in[0,\tau]}W_T(t)<-d\right\}
&\le\mathbb P\left\{\sup_{t\in[0,\tau]}(-X_T(t))>d/2\right\}\\
&\le\frac{2m_\tau}{d}.
\end{aligned}
\]
Here $\sup_{[0,\tau]}(-X_T)\ge0$ because $X_T(0)=0$.
Combining the estimates for the three regions and using reflection once more yields
\[
\limsup_{T\to\infty}\mathbb P\{L_T<-d\}
\le2\sum_{k=L}^{\infty}k^{-1-\eta}+\frac{4m_\tau}{d}.
\]
Letting $\tau\downarrow0$ and then $L\to\infty$ proves
\eqref{probconv2}, since $W_T(0)=W_T(T)=0$ and hence $L_T\le0$.

We now return to the conditioning identities, and in view of
$R_V(t,T)=\sigma_V^2(t)+\kappa_V(t,T)\ge-a_T$ together with
$R_V(T-t,T)=S_T^2-R_V(t,T)$, we obtain
\[
-\frac{a_T}{S_T^2}\le q_T(t)\le1+\frac{a_T}{S_T^2},
\qquad 0\le t\le T.
\]
Put $m(x)=\min(0,\sqrt2x)$. Using the endpoint values together with this bound, we obtain
\[
L_T+m(x)-\frac{\sqrt2|x|a_T}{S_T^2}
\le M_T(x)\le m(x).
\]
In view of \eqref{probconv2}, we have $M_T(x)\to m(x)$ in probability
for every fixed $x$, and since $0\le e^{M_T(x)}\le1$, it also follows that
$\mathbb E[e^{M_T(x)}]\to e^{m(x)}$.
Moreover, we have 
\[
0\le\mathbb E[e^{M_T(x)}]h_T(x)
\le e^{m(x)-x/\sqrt2}=e^{-|x|/\sqrt2}.
\]
Dominated convergence in the first conditioning identity therefore yields
\[
\frac{\mathcal G_V(T)}{b_T}
\longrightarrow\int_{\mathbb R}e^{-|x|/\sqrt2}\,dx=2\sqrt2.
\]
In the light of $\Psi(S_T/\sqrt2)\sim\sqrt2b_T$ and
$\gamma_V(T)=S_T^2/2$, the preceding limit proves \eqref{con1:res}.

For $C<0$, if $x\le C/\sqrt2$, then $M_T(x)\le\sqrt2x\le C$,
whereas if $x>C/\sqrt2$, we have $m(x)>C$ and hence
$\mathbb P\{M_T(x)>C\}\to1$.
Since the integrand in the second conditioning identity is bounded by
$\mathbf1_{\{x>C/\sqrt2\}}e^{-x/\sqrt2}$, which is integrable,
another application of dominated convergence yields
\[
\frac{p_V(T,C)}{b_T}
\longrightarrow\int_{C/\sqrt2}^{\infty}e^{-x/\sqrt2}\,dx
=\sqrt2e^{-C/2}.
\]
Using the same calculation for the endpoint probability, we obtain
\[
\frac{\mathbb P\{\sqrt2V(T)-S_T^2>C\}}{b_T}
=\int_{C/\sqrt2}^{\infty}h_T(x)\,dx
\longrightarrow\sqrt2e^{-C/2}.
\]
This proves \eqref{eqpersistence}. \QED \\

\prooftheo{ThmKrzys}
The claim follows by an argument similar to that used in the proof of
\cite[Thm 3.1]{debicki2017approximation}.
The stationarity of increments of $V_i$, $i=1,2$, implies that its
variogram $\gamma_{V_i}$ is conditionally negative definite. Hence
Schoenberg's theorem implies that for each $u>0$
\[
 r_u^{(i)}(s,t):=\exp\left(-\frac{2}{u^2}\gamma_{V_i}(s-t)\right),
 \qquad s,t\in[0,T],\quad i=1,2
\]
is a positive definite function and thus a valid covariance function.
Fix a finite set $K\subset[0,T]$ containing $0$, and let
$(X_u^{(i)}(t))_{t\in K}$ be a centered Gaussian vector with covariance
$r_u^{(i)}$. Its coordinates have unit variance. In view of the evenness
of variograms and \eqref{gammas}
\[
 r_u^{(1)}(s,t)\ge r_u^{(2)}(s,t),\qquad s,t\in K.
\]
Slepian's inequality, see e.g., \cite[Thm 3]{Lif} and \cite[Lem 4]{PI23},
therefore yields for $u>0$
\begin{equation}\label{slepian1}
\begin{split}
 \mathbb P\left\{\min_{t\in K}X_u^{(1)}(t)>u\right\}
 &\ge \mathbb P\left\{\min_{t\in K}X_u^{(2)}(t)>u\right\},\\
 \mathbb P\left\{\max_{t\in K}X_u^{(1)}(t)>u\right\}
 &\le \mathbb P\left\{\max_{t\in K}X_u^{(2)}(t)>u\right\}.
\end{split}
\end{equation}

Set $Y_i(t)=V_i(t)-V_i(0)$ and
$W_i(t)=\sqrt2Y_i(t)-2\gamma_{V_i}(t)$. We next establish the finite-grid
high-threshold limits
\begin{equation}\label{from.uniform}
\begin{split}
 \lim_{u\to\infty}
 \frac{\mathbb P\{\min_{t\in K}X_u^{(i)}(t)>u\}}{\Psi(u)}
 &=\mathbb E\left\{\min_{t\in K}e^{W_i(t)}\right\},\\
 \lim_{u\to\infty}
 \frac{\mathbb P\{\max_{t\in K}X_u^{(i)}(t)>u\}}{\Psi(u)}
 &=\mathbb E\left\{\max_{t\in K}e^{W_i(t)}\right\}.
\end{split}
\end{equation}
For this purpose fix $i$, write $\gamma=\gamma_{V_i}$, and suppress
the process index. Set $a_u(t)=r_u(t,0)$, $b_u(t)=u^2(a_u(t)-1)$ and define the centered
Gaussian residual
\[
 Z_u(t)=u\bigl(X_u(t)-a_u(t)X_u(0)\bigr),\qquad t\in K.
\]
This residual is independent of $X_u(0)$ and
\[
\begin{split}
 \operatorname{Cov}(Z_u(s),Z_u(t))
 &=u^2\bigl(r_u(s,t)-r_u(s,0)r_u(t,0)\bigr)\\
 &\longrightarrow 2\bigl(\gamma(s)+\gamma(t)-\gamma(s-t)\bigr)
 =2\operatorname{Cov}(Y(s),Y(t)).
\end{split}
\]
Since $a_u(t)\to1$ and $b_u(t)\to-2\gamma(t)$, conditioning on
$X_u(0)=u+w/u$ yields
\[
\begin{split}
 \bigl(u(X_u(t)-u)\bigr)_{t\in K}
 &\stackrel{d}{=}
 \bigl(Z_u(t)+a_u(t)w+b_u(t)\bigr)_{t\in K}\\
 &\Longrightarrow \bigl(W(t)+w\bigr)_{t\in K}.
\end{split}
\]
Here the equality in distribution refers to the conditional law.
For $w\ne0$, every limiting coordinate has a continuous distribution
unless $\gamma(t)=0$, in which case it equals $w$. Thus the probabilities
of a strictly positive minimum or maximum converge for every $w\ne0$,
without any nondegeneracy assumption on the Gaussian vectors.

Let $\mathcal L(x)=\min_{t\in K}x_t$ or
$\mathcal L(x)=\max_{t\in K}x_t$, and let $\varphi$ denote the standard
normal density. Conditioning as above implies
\[
\begin{split}
 \frac{\mathbb P\{\mathcal L(X_u)>u\}}{\Psi(u)}
 &=\frac{\varphi(u)}{u\Psi(u)}
 \int_{\mathbb R}e^{-w-w^2/(2u^2)}\\
 &\quad{}\times\mathbb P\left\{
 \mathcal L\bigl((Z_u(t)+a_u(t)w+b_u(t))_{t\in K}\bigr)>0
 \right\}\,dw.
\end{split}
\]
The prefactor tends to $1$. On $w\ge0$, the integrand is bounded by
$e^{-w}$. For the minimum, it vanishes on $w\le0$, since $0\in K$ and
the coordinate at $0$ equals $w$.
For the maximum, put $M=\max_{t\in K}\gamma(t)$. If $M=0$, all coordinates
of $X_u$ coincide almost surely and $W(t)=0$ on $K$, so both limits in
\eqref{from.uniform} equal $1$. If $M>0$, then, for all sufficiently
large $u$ and all $t\in K$
\[
 a_u(t)\ge\frac12,\qquad b_u(t)\le0,\qquad
 \operatorname{Var}(Z_u(t))
 =u^2\bigl(1-e^{-4\gamma(t)/u^2}\bigr)\le4M.
\]
Consequently, for $w<0$ the probability in the integrand for the maximum
is at most
\[
 \sum_{t\in K}\mathbb P\{Z_u(t)>|w|/2\}
 \le |K|\exp\left(-\frac{w^2}{32M}\right).
\]
The resulting bound $|K|\exp(|w|-w^2/(32M))$ is integrable on
$(-\infty,0)$. Dominated convergence and Tonelli's theorem now imply
\[
 \lim_{u\to\infty}
 \frac{\mathbb P\{\mathcal L(X_u)>u\}}{\Psi(u)}
 =\int_{\mathbb R}e^{-w}\mathbb P\{\mathcal L(W)+w>0\}\,dw
 =\mathbb E e^{\mathcal L(W)},
\]
which proves both assertions in \eqref{from.uniform}.

Combining \eqref{slepian1} and \eqref{from.uniform} yields
\[
\begin{split}
 \mathbb E\left\{\min_{t\in K}e^{W_1(t)}\right\}
 &\ge \mathbb E\left\{\min_{t\in K}e^{W_2(t)}\right\},\\
 \mathbb E\left\{\max_{t\in K}e^{W_1(t)}\right\}
 &\le \mathbb E\left\{\max_{t\in K}e^{W_2(t)}\right\}.
\end{split}
\]
Finally, take nested finite grids containing $0,T$ with dense union in
$[0,T]$. Sample continuity, dominated convergence for the minima, which
are bounded by $e^{W_i(0)}=1$ and monotone convergence for the maxima
imply
\[
 \mathcal G_{Y_1}(T)\ge\mathcal G_{Y_2}(T),\qquad
 \mathcal H_{Y_1}(T)\le\mathcal H_{Y_2}(T).
\]
In view of \Cref{PegAs}, $\mathcal G_{Y_i}(T)=\mathcal G_{V_i}(T)$; the
same spectral shift identity, with the infimum replaced by the supremum,
implies $\mathcal H_{Y_i}(T)=\mathcal H_{V_i}(T)$. Hence \eqref{nachsitz}
follows.
\QED
\\

\kE{
\proofkorr{cor:BM-dominance-G} We first prove (i). Setting \(Y(t)=Z(t)-Z(0)\), \(t\ge0\) we have  
$\gamma_Y=\gamma_Z$ and consequently \Cref{PegAs} implies
\[
\mathcal G_Y(T)=\mathcal G_Z(T),\qquad T>0.
\]
Since \(B\) and \(Y\) are independent we obtain
\[
\mathcal G_V(T)=\mathcal G_{B+Y}(T)
\ge
\mathcal G_B(T)\mathcal G_Y(T)
=
\mathcal G_B(T)\mathcal G_Z(T).
\] Since further
\[
\gamma_{B+Y}(t)=\gamma_B(t)+\gamma_Y(t)\ge \gamma_B(t),\qquad t\ge0
\]
 \Cref{ThmKrzys} yields
\[
\mathcal G_V(T)=\mathcal G_{B+Y}(T)\le \mathcal G_B(T).
\]
In the light of \eqref{as_G_st}, we have
\[
\mathcal G_Z(T)
=
e^{-2\sigma\sqrt{\ln T}+O(1)},\qquad T\to\infty,
\]
\aa{which leads to} \eqref{eqH}.\\
We next prove (ii). Put
\[
q_c(t):=\sqrt{t^2+c\sin^2(t/2)}=2\gamma_{V_c}(t),\qquad t\ge0.
\]
Since $q_c(t)\ge t$ we have
\[
\gamma_{V_c}(t)\ge \gamma_B(t)=\frac t2,\qquad t\ge0,
\]
and therefore \Cref{ThmKrzys} implies
\[
\mathcal G_{V_c}(T)\le \mathcal G_B(T).
\]
For the lower bound write
\[
q_c(t)=t+d_c(t),\qquad d_c(t):=q_c(t)-t.
\]
Then \(d_c(0)=0\), and for \(t>0\),
\[
d_c(t)
=
\frac{c\sin^2(t/2)}{q_c(t)+t}.
\]
Hence there exists \(\sigma\in(0,\infty)\) such that
\[
0\le \frac12d_c(t)\le \sigma^2\min(t,1),\qquad t\ge0.
\]
Let $W$ be a standard Brownian motion independent of $B$, and let
$S(t)=W(t+1)-W(t)$, $t\ge0$ be the corresponding Slepian process. Setting 
\[
Y(t)=\sigma(S(t)-S(0))
\]
then
\[
\gamma_Y(t)=\sigma^2\min(t,1),
\]
and consequently
\[
\gamma_{V_c}(t)
=
\frac t2+\frac12d_c(t)
\le
\gamma_B(t)+\gamma_Y(t)
=
\gamma_{B+Y}(t),\qquad t\ge0.
\]
The \aa{rest of the} proof follows as in the part (i), and is therefore omitted.\\
We continue with the proof of (iii).  Put 
\[
D(t):=\int_0^t sR(s)\,ds+t\int_t^\infty R(s)\,ds,\qquad t\ge0.
\]
Since
\[
\gamma_V(t)=\int_0^t(t-s)R(s)\,ds,
\]
we have
\[
\gamma_V(t)=at-D(t).
\]
Moreover, by the definition
\[
g(t)=\int_t^\infty(s-t)R(s)\,ds,
\]
we get
\[
m-g(t)=D(t),\qquad t\ge0.
\]
Since $t\mapsto g(|t|)$ is positive definite by assumption, there exists
a centered stationary Gaussian process $Y$, independent of $V$, with
covariance
\[
\Cov(Y(0),Y(t))=g(|t|),\qquad t\in\mathbb R.
\]
Since \(g(0)=m\), its variogram is
\[
\gamma_Y(t)=g(0)-g(t)=m-g(t)=D(t).
\]
Consequently, we have 
\[
\gamma_{V+Y}(t)=\gamma_V(t)+\gamma_Y(t)=at.
\]
If \(B\) is standard Brownian motion, then \(\gamma_{B(2a\cdot)}(t)=at\).
In view of \Cref{PegAs}, this implies
\[
\mathcal G_{V+Y}(T)=\mathcal G_B(2aT).
\]
Since \(D(t)\ge0\), we have
\[
\gamma_V(t)\le at=\gamma_{B(2a\cdot)}(t).
\]
An application of \Cref{ThmKrzys} therefore yields
\[
\mathcal G_V(T)\ge \mathcal G_B(2aT).
\]
Since the independence of $Y$ and $V$ implies
\[
\mathcal G_{V+Y}(T)\ge \mathcal G_V(T)\mathcal G_Y(T), \qquad T>0
\]
and thus
\[
\mathcal G_V(T)
\le
\frac{\mathcal G_B(2aT)}{\mathcal G_Y(T)}
\]
it remains only to estimate \(\mathcal G_Y(T)\).  Note that
\[
g(0)-g(t)=D(t)\sim at,\qquad t\downarrow0.
\]
Moreover, $D(h)\le ah$ for $h\ge0$, and hence Gaussian fourth moments imply
\[
\mathbb E|Y(t)-Y(s)|^4=12D(|t-s|)^2\le12a^2|t-s|^2.
\]
Consequently, $Y$ admits a continuous version, while the assumption
$\ln(t)g(t)\to0$ implies its ergodicity; in the light of the local expansion
above, the Pickands condition holds with $A=a$ and $\alpha=1$, so an
application of \Cref{PropKE} yields 
\[
\mathcal G_Y(T)
\sim e^{-2\sqrt m\,\sqrt{\ln T}-m}, \qquad T\to\infty
\]
establishing the proof. 
\QED \\
 } 
\proofkorr{cor.comp}
We establish the assertion for $\mathcal G_{B_H}(T)$, since the
corresponding argument for $\mathcal H_{B_H}(T)$ follows along the same
lines.

Let $0<H_1<H<H_2\le 1$ and
\[
\lambda=\frac{H_2-H}{H_2-H_1}>0.
\]
Define next 
\[
V(t)=\lambda^{1/2} T^{H-H_1}B_{H_1}(t)+(1-\lambda)^{1/2}T^{H-H_2}B_{H_2}(t),
\qquad t\in[0,T].
\]
In order to apply \Cref{ThmKrzys} to the processes $B_H$ and $V$, it is enough to prove that
\begin{equation}\label{comp1}
\gamma_{B_H}(t)\le \gamma_V(t), \qquad t\in[0,T].
\end{equation}
Since
$
\gamma_{B_H}(t)=t^{2H}\ks{/2},t\ge 0$
and
\[
\gamma_V(t)=\left( \lambda T^{2(H-H_1)}t^{2H_1}+(1-\lambda)T^{2(H-H_2)}t^{2H_2} \right)\ks{/2},
\]
\eqref{comp1} is equivalent, for $t>0$, to
\[
t^{2(H-H_1)}
\le
\lambda T^{2(H-H_1)}+(1-\lambda)T^{2(H-H_2)}t^{2(H_2-H_1)}.
\]
Thus it suffices to show that
\[
f(t):=\lambda T^{2(H-H_1)}+(1-\lambda)T^{2(H-H_2)}t^{2(H_2-H_1)}-t^{2(H-H_1)}
\ge 0
\]
for all $t\in[0,T]$.

Now
\[
f(0)=\lambda T^{2(H-H_1)}>0,
\qquad
f(T)=0,
\]
and
\[
f'(t)
=
2(H-H_1)t^{2(H-H_1)-1}
\left[\left(\frac{t}{T}\right)^{2(H_2-H)}-1\right]
<0,
\qquad t\in(0,T)
\]
since $H_2>H$. Hence $f$ is decreasing on $[0,T]$, and therefore
$f(t)\ge 0$ for all $t\in[0,T]$. This proves \eqref{comp1}.
Consequently, \Cref{ThmKrzys} implies
\[
\mathcal{G}_{B_H}(T)\ge \mathcal{G}_{V}(T), \qquad T>0.
\]

Next, using
\[
\inf_{t\in[0,T]}(A_t+B_t)\ge \inf_{t\in[0,T]}A_t+\inf_{t\in[0,T]}B_t,
\]
we obtain
\begin{align*}
\inf_{t\in [0,T]} (\sqrt{2}V(t)-\sigma^2_V(t))
&\ge
\inf_{t\in [0,T]}
\bigl(\sqrt{2}\lambda^{1/2} T^{H-H_1}B_{H_1}(t)-\lambda T^{2(H-H_1)}t^{2H_1}\bigr)
\\
&\quad+
\inf_{t\in [0,T]}
\bigl(\sqrt{2}(1-\lambda)^{1/2}T^{H-H_2}B_{H_2}(t)-(1-\lambda)T^{2(H-H_2)}t^{2H_2}\bigr)
\end{align*}
and hence
\begin{equation}\label{simp2}
 \mathcal{G}_{V}(T)
 \ge
 \mathcal{G}_{\lambda^{1/2} T^{H-H_1}B_{H_1}}(T)
 \mathcal{G}_{(1-\lambda)^{1/2}T^{H-H_2}B_{H_2}}(T).
\end{equation}

Finally, in view of the self-similarity of fBM, for $a>0$ and $H\in(0,1]$
we obtain
\[
\mathcal{G}_{aB_H}(T)=\mathcal{G}_{B_H}(a^{1/H}T).
\]
Combining this with \eqref{simp2} yields
\[
 \mathcal{G}_{V}(T)
 \ge
\mathcal{G}_{B_{H_1}}\left(\lambda^{1/(2H_1)}T^{H/H_1}\right)
\mathcal{G}_{B_{H_2}}\left((1-\lambda)^{1/(2H_2)}T^{H/H_2}\right),
\]
establishing the proof.
\QED

\prooftheo{ThmSlava2}
Fix \(\delta>0\) and set
\[
T_\delta:=\delta\left\lfloor \frac{T}{\delta}\right\rfloor .
\]
Since \(0,T_\delta\in[0,T_\delta]\cap\delta\mathbb Z\), we have
\[
\inf_{t\in[0,T_\delta]}
e^{\sqrt2V(t)-\sigma_V^2(t)}
\le
\inf_{t\in[0,T]\cap\delta\mathbb Z}
e^{\sqrt2V(t)-\sigma_V^2(t)}
\le
\min\left(1,e^{\sqrt2V(T_\delta)-\sigma_V^2(T_\delta)}\right).
\]
Taking expectations and using \eqref{boundBy}, we obtain
\[
\mathcal G_V(T_\delta)
\le
\mathcal G_V^\delta(T)
\le
2\Psi\left(\sqrt{\gamma_V(T_\delta)}\right).
\]
In view of \Cref{ThmSlava}, we have
\[
\mathcal G_V(T_\delta)
\sim
2\Psi\left(\sqrt{\gamma_V(T_\delta)}\right),
\qquad T\to\infty .
\]
Combining these relations, we obtain
\[
\mathcal G_V^\delta(T)
\sim
\mathcal G_V(T_\delta)
\sim
2\Psi\left(\sqrt{\gamma_V(T_\delta)}\right),
\qquad T\to\infty ,
\]
which proves \eqref{con12}.
Using the same endpoint comparison, we obtain
\[
p_V(T_\delta,C)
\le
\mathbb P\left\{
\inf_{t\in[0,T]\cap\delta\mathbb Z}
(\sqrt2V(t)-\sigma_V^2(t))>C
\right\}
\le
\mathbb P\left\{
\sqrt2V(T_\delta)-\sigma_V^2(T_\delta)>C
\right\}.
\]
In the light of \Cref{ThmSlava}, the lower and upper bounds satisfy
\[
p_V(T_\delta,C)
\sim
\mathbb P\left\{
\sqrt2V(T_\delta)-\sigma_V^2(T_\delta)>C
\right\},
\qquad T\to\infty .
\]
Consequently,
\[
\mathbb P\left\{
\inf_{t\in[0,T]\cap\delta\mathbb Z}
(\sqrt2V(t)-\sigma_V^2(t))>C
\right\}
\sim
p_V(T_\delta,C),
\qquad T\to\infty ,
\]
which proves \eqref{eqpersistence2}.

We now turn to the Brownian lattice case. Let
\[
G_N^\delta
:=
\mathcal G_{B_{1/2}}([0,\delta N]\cap\delta\mathbb Z).
\]
Let \(X\sim N(-\delta,2\delta)\), and let \(X_k\), \(k\ge1\), be iid copies of
\(X\). Setting
\[
S_0=0,\qquad S_k=\sum_{\ell=1}^kX_\ell,\quad k\ge1,
\]
we have
\[
G_N^\delta
=
\mathbb E\left\{e^{\min_{0\le k\le N}S_k}\right\}.
\]
Indeed,
\[
S_k\stackrel{d}{=}\sqrt2B(\delta k)-\delta k,
\qquad k\ge0.
\]

Let \(X'_k\), \(k\ge1\), be iid \(N(0,2\delta)\) random variables, and put
\[
S'_0=0,\qquad S'_k=\sum_{\ell=1}^kX'_\ell.
\]
The Gaussian change of measure, equivalently Girsanov's theorem, yields
\[
\mathbb E\left\{e^{\min_{0\le k\le N}S_k}\right\}
=
e^{-\delta N/4}
\mathbb E\left\{
e^{\min_{0\le k\le N}S'_k-S'_N/2}
\right\}.
\]
Let
$
X''_k=-X'_k, S''_0=0, S''_k=\sum_{\ell=1}^kX''_\ell=-S'_k,
$ and define
\[
R''_N:=\max_{0\le k\le N}S''_k.
\]
Then
\[
\mathbb E\left\{
e^{\min_{0\le k\le N}S'_k-S'_N/2}
\right\}
=
\mathbb E\left\{e^{-R''_N+S''_N/2}\right\}.
\]
Consequently,
\[
G_N^\delta
=
e^{-\delta N/4}A_N^\delta,
\qquad
A_N^\delta:=
\mathbb E\left\{e^{-R''_N+S''_N/2}\right\}.
\]

In view of Spitzer's formula, see for instance \cite{wendel1958spitzer}, for
real \(\alpha,\beta\) and \(|z|<1\),
\[
\sum_{n=0}^{\infty}
\mathbb E\left\{e^{i\alpha R''_n+i\beta S''_n}\right\}z^n
=
\exp\left(
\sum_{n=1}^{\infty}
\mathbb E\left\{
e^{i\alpha(S''_n\vee0)+i\beta S''_n}
\right\}
\frac{z^n}{n}
\right).
\]
We shall use this identity at \(\alpha=i\), \(\beta=-i/2\), with the
substitution justified coefficientwise. For each fixed \(N\), the
coefficient of \(z^N\) on the right-hand side is a finite polynomial in
\[
\mathbb E\left\{
 e^{i\alpha(S''_k\vee0)+i\beta S''_k}
\right\},\qquad 1\le k\le N.
\]
These expectations and the coefficient on the left-hand side are entire
in each of \(\alpha,\beta\), since the absolute values of the relevant
maxima and endpoints are bounded by
\(\sum_{k=1}^N|X''_k|\), which has exponential moments of every order.
For a fixed real \(\beta\), the one-variable identity theorem first
extends the coefficient identity to all complex \(\alpha\), and then,
for each such \(\alpha\), the same argument extends it to all complex
\(\beta\). At the required values,
\[
-R''_n+S''_n/2\le0,
\qquad -(S''_n\vee0)+S''_n/2=-|S''_n|/2\le0.
\]
Consequently, both resulting power series converge absolutely for
\(|z|<1\), and summing the coefficient identities yields
\begin{align}
\sum_{n=0}^{\infty}
\mathbb E\left\{e^{-R''_n+S''_n/2}\right\}z^n
&=
\exp\left(
\sum_{n=1}^{\infty}
\mathbb E\left\{
e^{-(S''_n\vee0)+S''_n/2}
\right\}
\frac{z^n}{n}
\right) \notag\\
&=
\exp\left(
\sum_{n=1}^{\infty}
\mathbb E\left\{
e^{-|S''_n|/2}
\right\}
\frac{z^n}{n}
\right).
\label{Spitzer}
\end{align}

Set
\[
f_\delta(n):=
\frac{1}{n}\mathbb E\left\{e^{-|S''_n|/2}\right\},
\qquad n\ge1.
\]
Since \(S''_n\sim N(0,2\delta n)\), completing the square implies
\begin{align*}
\mathbb E\left\{e^{-|S''_n|/2}\right\}
&=
\frac{2}{\sqrt{4\pi\delta n}}
\int_0^\infty
\exp\left(-\frac{x}{2}-\frac{x^2}{4\delta n}\right)\,dx\\
&=2e^{\delta n/4}\Psi\left(\sqrt{\delta n/2}\right)
\sim\frac{2}{\sqrt{\pi\delta n}},
\qquad n\to\infty .
\end{align*}
Hence
\begin{equation}
f_\delta(n)
\sim
\frac{2}{\sqrt{\pi\delta}}\frac{1}{n^{3/2}},
\qquad n\to\infty .
\label{f-asymp}
\end{equation}
In particular,
\[
M_\delta:=\sum_{n=1}^{\infty}f_\delta(n)<\infty.
\]
For the coefficient argument, write \(f=f_\delta\) and \(M=M_\delta\),
and choose \(0<b_1\le b_2<\infty\) such that
\[
b_1n^{-3/2}\le f(n)\le b_2n^{-3/2},\qquad n\ge1.
\]
Comparing coefficients in \eqref{Spitzer}, we obtain, for \(N\ge1\),
\[
A_N^\delta
=\sum_{j\ge1}\frac{f^{*j}(N)}{j!}
=
\sum_{j\ge1}
\sum_{\substack{c_1,\ldots,c_j\ge1\\ c_1+\cdots+c_j=N}}
\frac{f_\delta(c_1)\cdots f_\delta(c_j)}{j!},
\]
where \(f^{*j}\) denotes the \(j\)-fold convolution on the positive
integers. In each composition of \(N\) into \(j\) positive parts, a
largest part is at least \(N/j\); bounding its factor and summing the
other factors over all positive integers therefore implies
\begin{equation}\label{lattice-convolution-bound}
f^{*j}(N)\le b_2j^{5/2}M^{j-1}N^{-3/2},
\qquad j,N\ge1.
\end{equation}
The factor \(j\) accounts for the possible positions of the largest
part, and
\(\sum_{j\ge1}j^{5/2}M^{j-1}/j!<\infty\), so this estimate also
controls the sum over the number of parts.

Let \(B_{N,j}\) denote the sum of the products
\(f(c_1)\cdots f(c_j)\) over compositions with every part smaller than
\(3N/4\). For \(j=1\), this sum is zero. For a fixed \(j\ge2\),
a largest part is at least \(N/j\), whereas the sum of the remaining
parts is at least \(N/4\); consequently,
\[
B_{N,j}\le b_2j^{5/2}N^{-3/2}
\sum_{\substack{c_1,\ldots,c_{j-1}\ge1\\
 c_1+\cdots+c_{j-1}\ge N/4}}
f(c_1)\cdots f(c_{j-1})
=o(N^{-3/2}).
\]
Indeed, the unrestricted sum has the finite value \(M^{j-1}\), so its
tail tends to zero. In the light of
\eqref{lattice-convolution-bound}, dominated convergence in the sum over
\(j\) now yields
\[
\sum_{j\ge1}\frac{B_{N,j}}{j!}=o(N^{-3/2}).
\]
Every remaining composition has a unique part at least \(3N/4\).
Therefore
\begin{align*}
A_N^\delta
&=
\sum_{j\ge1}
\sum_{\substack{c_1,\ldots,c_j\ge1\\ c_1+\cdots+c_j=N\\
\exists \ell:\ c_\ell\ge 3N/4}}
\frac{f_\delta(c_1)\cdots f_\delta(c_j)}{j!}
+o(N^{-3/2})
\\
&=
\sum_{j\ge1}
\sum_{\substack{c_1,\ldots,c_j\ge1\\ c_1+\cdots+c_j=N\\
c_1\ge 3N/4}}
\frac{f_\delta(c_1)\cdots f_\delta(c_j)}{(j-1)!}
+o(N^{-3/2})
\\
&=
\sum_{j\ge1}
\sum_{\substack{c_2,\ldots,c_j\ge1\\ c_2+\cdots+c_j\le N/4}}
\frac{
f_\delta\!\left(N-\sum_{\ell=2}^j c_\ell\right)
f_\delta(c_2)\cdots f_\delta(c_j)}
{(j-1)!}
+o(N^{-3/2}).
\end{align*}
For each integer \(0\le s\le N/4\), the preceding two-sided bounds imply
\[
\frac{f_\delta(N-s)}{f_\delta(N)}
\le\frac{b_2}{b_1}\left(\frac43\right)^{3/2},
\]
whereas regular variation implies that this ratio tends to one for each
fixed nonnegative integer \(s\). Moreover,
\[
\sum_{j\ge1}\frac1{(j-1)!}
\sum_{c_2,\ldots,c_j\ge1}
 f_\delta(c_2)\cdots f_\delta(c_j)
=\sum_{j\ge1}\frac{M_\delta^{j-1}}{(j-1)!}
=e^{M_\delta}.
\]
Here and below, an empty product is interpreted as one and the sum of
an empty tuple as zero.
Thus one may first restrict \(c_2+\cdots+c_j\le K\), where the ratio
converges uniformly over the finitely many tuples, and then let
\(K\to\infty\), since the displayed bound and the finite total mass
control the omitted terms uniformly in \(N\). The last composition
expansion therefore yields
\[
A_N^\delta
=
f_\delta(N)
\sum_{j\ge0}
\sum_{c_1,\ldots,c_j\ge1}
\frac{f_\delta(c_1)\cdots f_\delta(c_j)}{j!}
+o(N^{-3/2})
=
f_\delta(N)
\exp\left(\sum_{k=1}^{\infty}f_\delta(k)\right)
+o(N^{-3/2}).
\]

Combining this with \eqref{f-asymp}, we get
\[
A_N^\delta
\sim
\frac{2}{\sqrt{\pi\delta}}\,
\frac{1}{N^{3/2}}
\exp\left(\sum_{k=1}^{\infty}f_\delta(k)\right),
\qquad N\to\infty .
\]
Therefore
\[
G_N^\delta
\sim
\frac{2e^{-\delta N/4}}{N^{3/2}\sqrt{\pi\delta}}\,
\exp\left(\sum_{k=1}^{\infty}f_\delta(k)\right),
\qquad N\to\infty .
\]
Writing
\[
F_\delta:=
\exp\left(\sum_{k=1}^{\infty}f_\delta(k)\right),
\]
this becomes
\[
G_N^\delta
\sim
\frac{2F_\delta e^{-\delta N/4}}{N^{3/2}\sqrt{\pi\delta}}
=
c_\delta\frac{e^{-\delta N/4}}{(\delta N)^{3/2}},
\qquad N\to\infty ,
\]
where
\[
c_\delta:=\frac{2\delta}{\sqrt\pi}F_\delta
=\frac{2\delta}{\sqrt\pi}
\exp\left\{2\sum_{k=1}^{\infty}
\frac{e^{\delta k/4}}{k}
\Psi\left(\sqrt{\delta k/2}\right)\right\}.
\]

It remains to identify \(F_\delta\) and to prove the limit of \(c_\delta\).
Letting \(z\uparrow1\) in \eqref{Spitzer} is justified by monotone convergence,
because all coefficients on the left-hand side are non-negative and because
\(\sum_{k=1}^{\infty}f_\delta(k)<\infty\). Thus
\[
F_\delta
=
\sum_{n=0}^{\infty}
\mathbb E\left\{e^{-R''_n+S''_n/2}\right\}
=
\sum_{n=0}^{\infty}A_n^\delta .
\]
Since
\[
G_n^\delta=e^{-\delta n/4}A_n^\delta,
\]
we obtain
\[
F_\delta
=
\sum_{n=0}^{\infty}
e^{\delta n/4}G_n^\delta
=
\sum_{n=0}^{\infty}
e^{\delta n/4}
\mathcal G_{B_{1/2}}([0,\delta n]\cap\delta\mathbb Z).
\]
Consequently,
\[
c_\delta
=
\frac{2\delta}{\sqrt\pi}
\sum_{n=0}^{\infty}
e^{\delta n/4}
\mathcal G_{B_{1/2}}([0,\delta n]\cap\delta\mathbb Z).
\]

We now show that
\[
\delta F_\delta\to4,
\qquad \delta\downarrow0.
\]
Indeed,
\begin{align}
\delta F_\delta
&=
\delta
\sum_{n=0}^{\infty}
e^{\delta n/4}
\mathcal G_{B_{1/2}}([0,\delta n]\cap\delta\mathbb Z)
\notag\\
&=
\sum_{n=0}^{\infty}
\int_{\delta n}^{\delta(n+1)}
e^{\delta n/4}
\mathcal G_{B_{1/2}}([0,\delta n]\cap\delta\mathbb Z)\,dx
\notag\\
&=
\int_0^\infty
e^{\delta\lfloor x/\delta\rfloor/4}
\mathcal G_{B_{1/2}}
\left([0,x]\cap\delta\mathbb Z\right)\,dx .
\label{delta-F-integral}
\end{align}
For each fixed \(x\ge0\), the continuity of Brownian paths implies
\[
\inf_{t\in[0,x]\cap\delta\mathbb Z}(\sqrt2B(t)-t)
\to
\inf_{t\in[0,x]}(\sqrt2B(t)-t),
\qquad \delta\downarrow0,
\]
almost surely. Since the corresponding exponentials are bounded by \(1\),
\[
\mathcal G_{B_{1/2}}([0,x]\cap\delta\mathbb Z)
\to
\mathcal G_{B_{1/2}}(x),
\qquad \delta\downarrow0.
\]
Also,
\[
\delta\left\lfloor \frac{x}{\delta}\right\rfloor\to x.
\]
Hence the integrand in \eqref{delta-F-integral} converges pointwise to
\[
e^{x/4}\mathcal G_{B_{1/2}}(x).
\]

It remains to justify domination. For \(0<\delta\le1\), let
\[
m_\delta:=\left\lceil\frac1\delta\right\rceil,
\qquad
\Delta_\delta:=m_\delta\delta\in[1,2].
\]
Then \(\Delta_\delta\mathbb Z\subset\delta\mathbb Z\). For
\[
n_\delta(x):=\left\lfloor\frac{x}{\delta}\right\rfloor,
\qquad
n_\Delta(x):=\left\lfloor\frac{x}{\Delta_\delta}\right\rfloor,
\]
we have
\[
\Delta_\delta n_\Delta(x)\le \delta n_\delta(x)\le x
\]
and therefore
\[
\mathcal G_{B_{1/2}}([0,x]\cap\delta\mathbb Z)
\le
\mathcal G_{B_{1/2}}
\left([0,\Delta_\delta n_\Delta(x)]\cap\Delta_\delta\mathbb Z\right).
\]
Moreover,
\[
\delta n_\delta(x)
\le
\Delta_\delta n_\Delta(x)+\Delta_\delta,
\]
and hence
\[
e^{\delta n_\delta(x)/4}
\le
e^{\Delta_\delta/4}e^{\Delta_\delta n_\Delta(x)/4}
\le
e^{1/2}e^{\Delta_\delta n_\Delta(x)/4}.
\]
The required coefficient bound is uniform for \(\Delta\in[1,2]\).
Indeed, the Gaussian integral above implies
\[
f_\Delta(n)\le\frac{2}{\sqrt{\pi\Delta}}n^{-3/2}
\le Cn^{-3/2},\qquad
M_\Delta\le C\sum_{n=1}^{\infty}n^{-3/2}=:M_*<\infty,
\]
where \(C=2/\sqrt\pi\) is independent of \(\Delta\). Applying the
largest-part estimate \eqref{lattice-convolution-bound} with these
uniform upper bounds yields, for \(n\ge1\),
\[
A_n^\Delta
\le Cn^{-3/2}
\sum_{j\ge1}\frac{j^{5/2}M_*^{j-1}}{j!}.
\]
Since \(A_0^\Delta=1\), there is consequently a constant
\(C_1<\infty\) such that, for all \(\Delta\in[1,2]\) and all
\(n\ge0\),
\[
e^{\Delta n/4}
\mathcal G_{B_{1/2}}([0,\Delta n]\cap\Delta\mathbb Z)
\le
\frac{C_1}{(1+n)^{3/2}}.
\]
Since \(x<\Delta_\delta(n_\Delta(x)+1)\le2(n_\Delta(x)+1)\),
we obtain, for all \(0<\delta\le1\) and \(x\ge0\),
\[
e^{\delta\lfloor x/\delta\rfloor/4}
\mathcal G_{B_{1/2}}([0,x]\cap\delta\mathbb Z)
\le
\frac{C_2}{(1+x)^{3/2}},
\]
with some constant \(C_2<\infty\). The right-hand side is integrable on
\([0,\infty)\). Therefore, by dominated convergence,
\[
\delta F_\delta
\to
\int_0^\infty e^{x/4}\mathcal G_{B_{1/2}}(x)\,dx,
\qquad \delta\downarrow0.
\]
In view of \eqref{BMcase}, integration up to \(T>0\) yields
\[
\int_0^T e^{x/4}\mathcal G_{B_{1/2}}(x)\,dx
=4+4(T-2)e^{T/4}\Psi\left(\sqrt{T/2}\right)
 -4\sqrt{T/\pi}.
\]
For example, this identity follows by differentiating the right-hand
side, whose value at \(T=0\) is zero. The normal tail expansion
\(e^{T/4}\Psi(\sqrt{T/2})
=(\pi T)^{-1/2}(1+O(T^{-1}))\) now implies
\[
\int_0^\infty e^{x/4}\mathcal G_{B_{1/2}}(x)\,dx=4.
\]
Hence
\[
\delta F_\delta\to4,
\qquad \delta\downarrow0.
\]
Finally,
\[
c_\delta=\frac{2\delta}{\sqrt\pi}F_\delta,
\]
and therefore
\[
c_\delta\to\frac{8}{\sqrt\pi},
\qquad \delta\downarrow0.
\]
This completes the proof.
\QED

\section*{Acknowledgments}  We thank two anonymous referees and the  Associate Editor for their valuable comments and suggestions.
K. D\c{e}bicki
was partially supported by the National Science Centre, Poland,  Grant
No 2024/55/B/ST1/01062
(2025-2028).

\bibliographystyle{ieeetr}

\bibliography{GG_V15}

\begin{thebibliography}{10}

\bibitem{kab2009}
Z.~Kabluchko, M.~Schlather, and L.~de~Haan, ``Stationary max-stable fields associated to negative definite functions,'' {\em Ann. Probab.}, vol.~37, pp.~2042--2065, 2009.

\bibitem{kulik:soulier:2020}
R.~Kulik and P.~Soulier, {\em Heavy tailed time series}.
\newblock Springer, Cham, 2020.

\bibitem{Hrovje}
H.~Planini\'{c} and P.~Soulier, ``The tail process revisited,'' {\em Extremes}, vol.~21, no.~4, pp.~551--579, 2018.

\bibitem{bro1977}
B.~M. Brown and S.~I. Resnick, ``Extreme values of independent stochastic processes,'' {\em J. Appl. Probab.}, vol.~14, pp.~732--739, 1977.

\bibitem{eddy1980distribution}
W.~F. Eddy, ``The distribution of the convex hull of a {G}aussian sample,'' {\em J. Appl. Probab.}, vol.~17, no.~3, pp.~686--695, 1980.

\bibitem{Resnickart}
S.~Resnick, {\em The Art of Finding Hidden Risks: Hidden Regular Variation in the 21st Century}.
\newblock Springer Nature, 2024.

\bibitem{Ilya25}
B.~Basrak, N.~Milinčević, and I.~Molchanov, ``Foundations of regular variation on topological spaces,'' {\em ArXiv preprint 2503.00921}, 2025.

\bibitem{MR3745388}
K.~D{\c{e}}bicki and E.~Hashorva, ``On extremal index of max-stable stationary processes,'' {\em Probab. Math. Statist.}, vol.~37, no.~2, pp.~299--317, 2017.

\bibitem{debicki2017approximation}
K.~D{\c{e}}bicki and E.~Hashorva, ``Approximation of supremum of max-stable stationary processes \& {P}ickands constants,'' {\em J. Theoret. Probab.}, vol.~33, no.~1, pp.~444--464, 2020.

\bibitem{ZKE}
K.~D{\c{e}}bicki, E.~Hashorva, and Z.~Michna, ``On the continuity of {P}ickands constants,'' {\em J. Appl. Probab.}, vol.~59, no.~1, pp.~187--201, 2022.

\bibitem{Htilt}
E.~Hashorva, ``Representations of max-stable processes via exponential tilting,'' {\em Stochastic Process. Appl.}, vol.~128, no.~9, pp.~2952--2978, 2018.

\bibitem{Genna04}
G.~Samorodnitsky, ``Extreme value theory, ergodic theory and the boundary between short memory and long memory for stationary stable processes,'' {\em Ann. Probab.}, vol.~32, no.~2, pp.~1438--1468, 2004.

\bibitem{Genna04c}
G.~Samorodnitsky, ``Maxima of continuous-time stationary stable processes,'' {\em Adv. in Appl. Probab.}, vol.~36, no.~3, pp.~805--823, 2004.

\bibitem{WangStoev}
Y.~Wang and S.~A. Stoev, ``On the structure and representations of max-stable processes,'' {\em Adv. in Appl. Probab.}, vol.~42, no.~3, pp.~855--877, 2010.

\bibitem{hashorva2021shiftinvariant}
E.~Hashorva, ``Shift-invariant homogeneous classes of random fields,'' {\em Journal of Mathematical Analysis and Applications}, p.~128517, 2024.

\bibitem{hashorva2025cluster}
E.~Hashorva, ``Cluster random fields and random-shift representations,'' {\em Journal of Theoretical Probability}, vol.~38, no.~3, p.~50, 2025.

\bibitem{DHL}
K.~D{\c{e}}bicki, E.~Hashorva, and P.~Liu, ``Uniform tail approximation of homogenous functionals of {G}aussian fields,'' {\em Advances in Applied Probability}, vol.~49, no.~4, pp.~1037--1066, 2017.

\bibitem{DeK14}
K.~D{\c{e}}bicki and K.~M. Kosi{\'n}ski, ``On the infimum attained by the reflected fractional {B}rownian motion,'' {\em Extremes}, vol.~17, no.~3, pp.~431--446, 2014.

\bibitem{dkebicki2016extremes}
K.~D{\c{e}}bicki and P.~Liu, ``Extremes of stationary gaussian storage models,'' {\em Extremes}, vol.~19, no.~2, pp.~273--302, 2016.

\bibitem{Molchan1999}
G.~Molchan, ``Maximum of fractional {B}rownian motion: probabilities of small values,'' {\em Communications in Mathematical Physics}, vol.~205, no.~1, pp.~97--111, 1999.

\bibitem{aurzada2013persistence}
F.~Aurzada and C.~Baumgarten, ``Persistence of fractional brownian motion with moving boundaries and applications,'' {\em Journal of Physics A: Mathematical and Theoretical}, vol.~46, no.~12, p.~125007, 2013.

\bibitem{AuS15}
F.~Aurzada and T.~Simon, ``Persistence probabilities and exponents,'' {\em L{\'e}vy matters V}, vol.~2149, pp.~183--224, 2015.

\bibitem{perB}
G.~Molchan, ``The persistence exponents of {G}aussian random fields connected by the {L}amperti transform,'' {\em Journal of Statistical Physics}, vol.~186, no.~2, p.~21, 2022.

\bibitem{perC}
F.~Aurzada and S.~Mukherjee, ``Persistence probabilities of weighted sums of stationary gaussian sequences,'' {\em Stochastic Processes and their Applications}, vol.~159, pp.~286--319, 2023.

\bibitem{FFN21}
N.~Feldheim, O.~Feldheim, and S.~Nitzan, ``Persistence of gaussian stationary processes,'' {\em The Annals of Probability}, vol.~49, no.~3, pp.~1067--1096, 2021.

\bibitem{aurzada2026}
F.~Aurzada and S.~Müller, ``Persistence probability of fractional {B}rownian motion with random {H}urst exponent,'' {\em ArXiv preprint arXiv:2603.14934}, 2026.

\bibitem{feldheim2025persistence}
N.~D. Feldheim, O.~N. Feldheim, and S.~Mukherjee, ``Persistence and ball exponents for {G}aussian stationary processes,'' {\em Communications on Pure and Applied Mathematics}, vol.~78, no.~10, pp.~1949--2000, 2025.

\bibitem{LiS01}
W.~V. Li and Q.-M. Shao, ``Lower tail probabilities for gaussian processes,'' {\em The Annals of Probability}, vol.~32, no.~1A, pp.~216--242, 2004.

\bibitem{DeS17}
A.~Dembo and S.~Mukherjee, ``Persistence of gaussian processes: non-summable correlations,'' {\em Probability Theory and Related Fields}, vol.~169, no.~3, pp.~1007--1039, 2017.

\bibitem{aurzada2018persistence}
F.~Aurzada, N.~Guillotin-Plantard, and F.~P{\`e}ne, ``Persistence probabilities for stationary increment processes,'' {\em Stochastic Processes and their Applications}, vol.~128, no.~5, pp.~1750--1771, 2018.

\bibitem{PickandsB}
J.~Pickands, III, ``Asymptotic properties of the maximum in a stationary {G}aussian process,'' {\em Trans. Amer. Math. Soc.}, vol.~145, pp.~75--86, 1969.

\bibitem{kabluchko2014limiting}
Z.~Kabluchko and Y.~Wang, ``Limiting distribution for the maximal standardized increment of a random walk,'' {\em Stochastic Processes and their Applications}, vol.~124, no.~9, pp.~2824--2867, 2014.

\bibitem{bisewski2025speed}
K.~Bisewski and G.~Jasnovidov, ``On the speed of convergence of discrete pickands constants to continuous ones,'' {\em Journal of Applied Probability}, vol.~62, no.~1, pp.~111--135, 2025.

\bibitem{BD}
G.~Baxter and M.~D. Donsker, ``On the distribution of the supremum functional for processes with stationary independent increments,'' {\em Transactions of the American Mathematical Society}, vol.~85, no.~1, pp.~73--87, 1957.

\bibitem{Led}
M.~R. Leadbetter, G.~Lindgren, and H.~Rootz{\'e}n, {\em Extremes and related properties of random sequences and processes}.
\newblock Springer Science \& Business Media, 2012.

\bibitem{Lif}
M.~A. Lifshits, {\em Gaussian random functions}, vol.~322.
\newblock Springer Science \& Business Media, 2013.

\bibitem{PI23}
P.~Ievlev, ``Extremes of locally-homogenous vector-valued {G}aussian processes,'' {\em Extremes}, pp.~1--27, 2024.

\bibitem{wendel1958spitzer}
J.~G. Wendel, ``Spitzer's formula: a short proof,'' {\em Proceedings of the American Mathematical Society}, vol.~9, no.~6, pp.~905--908, 1958.

\end{thebibliography}

\end{document}